\documentclass[11pt]{amsart}
\usepackage{amsmath,amsfonts,amssymb,amscd,amsthm}
\usepackage{mathrsfs,bbm}
\usepackage{stackengine}
\usepackage{extpfeil}
\usepackage{indentfirst}
\usepackage{latexsym}
\usepackage{enumerate}
\usepackage{float}
\usepackage{graphicx}
\usepackage{xcolor}
\usepackage{comment}
\usepackage{centernot}
\usepackage{bm}

\title[Success or Failure of Finite-Time Tracking in the Pursuer--Evader Model] {Finite-Time Tracking Success or Failure in the Pursuer--Evader Model}

\author[W. Shim]{Woojoo Shim}
\address[Woojoo Shim]{\newline Department of Mathematics Education \newline Kyungpook National University, Daegu 41566, Republic of Korea}
\email{wjshim@knu.ac.kr}

\author[H. Ahn]{Hyunjin Ahn$^*$}
\address[Hyunjin Ahn]{\newline Department of Data Technology \newline Myongji University, Seoul 03674, Republic of Korea}
\email{ahj92@mju.ac.kr}

\newtheorem{theorem}{Theorem}[section]
\newtheorem{lemma}{Lemma}[section]
\newtheorem{corollary}{Corollary}[section]
\newtheorem{proposition}{Proposition}[section]

\newtheorem{remark}{Remark}[section]

\newtheorem{definition}{Definition}[section]

\newcommand{\bbr}{\mathbb R}

\def\charf {\mbox{{\text 1}\kern-.30em {\text l}}}

\calclayout

\begin{document}
	\date{\today}
	
	
	\thanks{Acknowledgment: The work of H. Ahn was supported by the National Research Foundation of Korea (NRF) grant funded by the Korea government(MSIT) (2022R1C12007321).}
	
	\thanks{$^*$ Corresponding author}
	
\begin{abstract}
We investigate a sufficient framework for determining the success or failure of finite-time tracking in pursuer--evader dynamics. To this end, we design a nonsmooth coupled multi-agent system by adopting a previously established pursuer--evader model from the existing literature, which characterizes asymptotic tracking success and failure, as the backbone. We then provide appropriate sufficient conditions, expressed in terms of initial data, system parameters, and kernel functionals, under which some pursuers can successfully track evaders within a finite time, while others fail to do so. In particular, we show that the pursuer--evader subsystem associated with the pursuers that successfully track the evader, together with the corresponding evaders, exhibits finite-time stability.
\end{abstract}
	
	\maketitle
	
	\centerline{\date}
	

\vspace{-1cm}
	
%
%
%
	\section{Introduction} \label{sec:1}
	\setcounter{equation}{0}
Tracking of active or passive targets is a basic behavior widely observed in both nature and human activities. An active target tries to escape or avoid the trackers, while a passive target moves on its own without trying to hide or resist being tracked. For example, in nature, a lion chasing a gazelle or a cat hunting a mouse involves an active target because the prey tries to run away. In contrast, when a solar panel follows the movement of the sun, the sun is a passive target because it does not try to avoid being tracked. In military settings, a drone that follows an enemy aircraft performing evasive maneuvers represents active target tracking, while a missile moving toward a stationary radar station shows passive target tracking. In daily life, an automatic camera that follows a soccer player during a game tracks an active target, while a GPS guided delivery robot moving toward a customer’s home or a light that follows an airplane during landing tracks a passive target. Similarly, a self driving car adjusting its route toward a fixed destination shows passive target tracking, whereas a police drone chasing a suspect’s car that keeps changing direction is another example of active tracking.

For these reasons, the study of tracking active or passive targets has been a dynamic area of research within the mathematics and engineering communities. Examples include the asymptotic tracking of a moving or stationary target and collision avoidance on Riemannian manifolds \cite{A-B-H-Y}; pursuer--evader games \cite{D-A-G, Isa, L-C, S-E-K-S, W-C-C-H-C, W-P-G}; scenarios involving multiple pursuers \cite{F-W-X-C, V-C-G-M}; leader-follower formation control problems where the leaders move with a constant velocity \cite{T-T-V-N-A}; vehicle platoon problems \cite{P-M-M, S-S-M}; alternating tracking of a moving or stationary target on Riemannian manifolds \cite{A-H-1}; stochastic asymptotic tracking of a moving or stationary target \cite{A-H3}; finite-time consensus algorithms for multi-agent systems with integrator dynamics \cite{P-P-S}; design of asymptotic tracking controllers \cite{F}; hydrodynamic descriptions for target tracking \cite{A-H5}; active target defense based on missile dynamics \cite{G-C-P}; asymptotic tracking of infinitely moving or stationary targets \cite{A-H-2}; fixed pattern formation with collision avoidance \cite{Do}; stochastic descriptions of fixed pattern formations \cite{C-O-T}; and the emergence of partial asymptotic tracking in the pursuer--evader model \cite{A-B-H-Y1}.

However, the previous studies have not fully explored whether tracking an active target can succeed or fail within a finite time. In many real-world situations, the key question is not whether tracking will eventually occur, but whether it can be achieved within a limited time. Studying such finite-time tracking dynamics helps us understand how fast and efficiently a system can respond, offering several important advantages in control and stability. Compared with asymptotic approaches, finite-time dynamics often provide faster convergence, stronger robustness against disturbances, and better transient performance. Although the convergence time may depend on the initial conditions, the system still guarantees convergence within a finite period, which is highly valuable in time-critical or safety-sensitive applications. For instance, in security and defense systems, tracking and capturing a fleeing intruder within a finite time is essential for mission success. In robotics, finite-time control enables autonomous robots to detect and avoid obstacles or follow moving targets quickly and stably. In drone operations, real-time tracking of fast-moving delivery objects or surveillance subjects relies on such rapid and robust convergence. Similarly, in autonomous driving, the ability to track pedestrians or surrounding vehicles in real time ensures safety and prevents accidents. Therefore, analyzing the success or failure of finite-time tracking is not only of theoretical interest but also crucial for enhancing the responsiveness, robustness, and overall effectiveness of control systems in practical applications.

\vspace{0.12cm}

Therefore, we are primarily concerned with the following issue in this study:

\vspace{0.12cm}

\begin{quote}
	``How can an effective pursuer--evader model be designed to represent the success or failure of a pursuer in tracking an evader within a finite time? Furthermore, can we propose a rigorous framework to describe when the pursuer successfully tracks the evader or fails within a finite time?"
\end{quote}

To address this issue, we draw inspiration from the pursuer--evader model studied in \cite{A-B-H-Y1}, which characterizes the success or failure of pursuers in asymptotically tracking evaders. To introduce this model, we first set up the mechanical states of the pursuer agents $(x_i, v_i)$ and the evader agents $(x_i^*, v_i^*)$. The authors of \cite{A-B-H-Y1} proposed the following coupled multi-agent system for $\{(x_i, v_i)\}_{i=1}^N$ and $\{(x_i^*, v_i^*)\}_{i=1}^N$ representing the pursuer--evader dynamics:

\vspace{-0.34cm}

\begin{equation}\label{wolf}
	\begin{cases} 
		\vspace{0.2cm}
		\displaystyle {\dot{x}_i} =v_i,\quad t>0,\quad  i\in [N]:=\{1,\ldots,N\},\\
		\vspace{0.2cm}
		\displaystyle \dot{v}_i=\underbrace{\frac{\kappa_{f}}{N}\sum_{j=1}^{N}\phi_f(\|(x_i-x_j)-(x_i^*-x_j^*)\|)\left(v_j- v_i\right)}_{\text{flocking force for pursuers}}  \\
		\vspace{0.2cm}
		\displaystyle \qquad\quad \underbrace{-\gamma_i (v_i-v_i^*)}_{\text{velocity-tracking force}} \underbrace{- \kappa_{r,i} \phi_{r,i} (\|x_i-x_i^*\|^2 )   (x_i-x^*_i)}_{\text{position-attraction force}},\\
		\displaystyle (x_i, v_i)\bigg|_{t=0+} =(x_i(0),v_i(0))\in \mathbb{R}^d\times \mathbb{R}^d,
	\end{cases}
\end{equation}
\begin{equation}\label{rabbit}
	\begin{cases} 
		\vspace{0.2cm}
		\displaystyle {\dot{x}^*_i} =v^*_i,\quad t>0,\quad  i\in [N]:=\{1,\ldots,N\},\\
		\vspace{0.2cm}
		\displaystyle \dot{v}^*_i=\underbrace{\frac{\kappa_{f}}{N}\sum_{j=1}^{N}\phi_f(\|(x_i-x_j)-(x_i^*-x_j^*)\|)\left(v^*_j- v^*_i\right)}_{\text{flocking force for evaders}}\\
		\vspace{0.2cm}
		\displaystyle \qquad\quad \underbrace{-\gamma^*_i (v_i-v_i^*)}_{\text{velocity-evading force}}  \underbrace{- \kappa_{c,i} \phi_{c,i} (\|x_i-x_i^*\|^2 )   (x_i-x^*_i)}_{\text{position-repulsion force}},\\
		\displaystyle (x^*_i, v^*_i)\bigg|_{t=0+} =(x^*_i(0),v^*_i(0))\in \mathbb{R}^d\times \mathbb{R}^d,
	\end{cases}
\end{equation} 
where $N$ denotes the number of agents and $\|\cdot\|$ denotes the standard Euclidean $l_2$-norm. Next, $\kappa_f$ represents the strength of flocking force for pursuers and evaders, respectively. $\gamma_i$ denotes the strength of velocity-tracking force related to the $i$-th pursuer pursuing the $i$-th evader, while $\gamma^*_i$ denotes the strength of velocity-evading force related to the $i$-th evader escaping from the $i$-th pursuer. In what follows, $\kappa_{r,i}$ characterizes the strength of position-attraction force related to the $i$-th pursuer tracking the $i$-th evader, while $\kappa_{c,i}$ characterizes the strength of position-repulsion force related to the $i$-th evader evading the $i$-th pursuer. Further, $\phi_f$ denotes a kernel functional regarding the flocking force for pursuers and evaders, respectively, $\phi_{r,i}$ denotes a kernel functional that represents the position-attraction force exerted by the $i$-th pursuer to track the $i$-th evader, and $\phi_{c,i}$ denotes a kernel functional that represents the position-repulsion force exerted by the $i$-th evader to evade the $i$-th pursuer.

Here, the authors of \cite{A-B-H-Y1} provided a sufficient framework for determining the success or failure of asymptotic tracking between pursuers and evaders in the coupled system $\eqref{wolf}-\eqref{rabbit}$. For further details on previous results, we refer to Section \ref{sec:2.3}.

Now, motivated by $\eqref{wolf}-\eqref{rabbit}$, we seek to propose our main model. Initially, let $u$ be a vector and $p, \alpha$ be two real numbers satisfying 
\begin{align*} 
	p \geq 1, \quad 0 < \alpha < 1, \quad u = (u_1, \ldots, u_d) \in \bbr^d, 
\end{align*}
and define the $\ell^p$-norm of $u$, a sign function, and $\overrightarrow{\textup{sig}}_\alpha(u)$ as follows:
\begin{align}
	\begin{aligned} \label{A-0}
		&  \|u\|_{p}:=\left(\sum_{k=1}^d |u_k|^p \right)^{\frac{1}{p}},\quad \|u\|:=\|u\|_2,\quad
		\textup{sgn}(x):=
		\begin{cases}
			-1,\qquad & x<0, \\
			0,\qquad  & x=0,\\
			1,\qquad  & x>0, 
		\end{cases} \\
		& \overrightarrow{\textup{sig}}_\alpha(u):=(\textup{sgn}(u_1)|u_1|^\alpha,\ldots,~\textup{sgn}(u_d)|u_d|^\alpha).
	\end{aligned}
\end{align}
To explain the idea behind the derivation of the proposed coupled system to be introduced later, we consider the following Cauchy problem for an ordinary differential equation (ODE) with a vector field that is not locally Lipschitz continuous: For $\alpha \in (0,1)$,
\begin{equation*}
	\begin{cases}
		\displaystyle {\dot x} = -\mbox{sgn}(x) |x|^{\alpha}, \quad t > 0, \\
		\displaystyle x(0) = x^0 \neq 0.
	\end{cases}
\end{equation*}
Then, the explicit solution is given by the following form:
\[ x(t) = 
\begin{cases}
	\displaystyle \mbox{sgn}(x^0) \Big(  |x^0|^{1-\alpha} - t(1-\alpha)  \Big)^{\frac{1}{1-\alpha}}, \quad &0 < t \leq \tau^*, \\
	\displaystyle 0, \quad &t \geq \tau^*,
\end{cases}
\]
where $\tau^* = \frac{|x^0|^{1-\alpha}}{1-\alpha}$. Hence, the solution $ x = x(t) $ reaches zero at the finite time $ \tau^* $.
\[    \lim_{t \to \tau^*-} x(t) = 0. \]

Subsequently, based on the above motivating example, combining $\eqref{wolf}-\eqref{rabbit}$ with \eqref{A-0}, we propose the following nonsmooth coupled system that represents pursuer--evader interactions in terms of $\{(x_i, v_i)\}_{i=1}^N$ and $\{(x_i^*, v_i^*)\}_{i=1}^N$:
\begin{equation}\label{Fwolf}
	\begin{cases} 
		\vspace{0.2cm}
		\displaystyle {\dot{x}_i} =v_i,\quad t>0,\quad  i\in [N]:=\{1,\ldots,N\},\\
		\vspace{0.2cm}
		\displaystyle \dot{v}_i=\underbrace{\frac{\kappa_{f}}{N}\sum_{j=1}^{N}\phi_f(\|(x_i-x_j)-(x_i^*-x_j^*)\|_{\alpha_3+1})\overrightarrow{\textup{sig}}_{\alpha_1}((v_j-v_j^*)- (v_i-v_i^*))}_{\text{{{relative}} flocking force between pursuers and evaders}}  \\
		\vspace{0.2cm}
		\displaystyle\qquad\quad+\underbrace{\frac{\kappa_{f}}{N}\sum_{j=1}^{N}\phi_f(\|(x_i-x_j)-(x_i^*-x_j^*)\|_{\alpha_3+1})\overrightarrow{\textup{sig}}_{\alpha_1}(v_j^*-v_i^*)}_{\text{{{disrupting}} force for the flocking of evaders}}  \\
		\vspace{0.2cm}
		\displaystyle \qquad\quad \underbrace{-\gamma_i~\overrightarrow{\textup{sig}}_{\alpha_2} (v_i-v_i^*)}_{\text{velocity-tracking force}}\hspace{0.2cm} \underbrace{- \kappa_{r,i} \phi_{r,i} (\|x_i-x_i^*\|_{\alpha_3+1}^{\alpha_3+1} )   \overrightarrow{\textup{sig}}_{\alpha_3}(x_i-x^*_i)}_{\text{position-attraction force}},\\
		\displaystyle (x_i, v_i)\bigg|_{t=0+} =(x_i(0),v_i(0))\in \mathbb{R}^d\times \mathbb{R}^d,
	\end{cases}
\end{equation}
\begin{equation}\label{Frabbit}
	\begin{cases} 
		\vspace{0.2cm}
		\displaystyle {\dot{x}^*_i} =v^*_i,\quad t>0,\quad  i\in [N]:=\{1,\ldots,N\},\\
		\vspace{0.2cm}
		\displaystyle \dot{v}^*_i=\underbrace{\frac{\kappa_{f}}{N}\sum_{j=1}^{N}\phi_f(\|(x_i-x_j)-(x_i^*-x_j^*)\|_{\alpha_3+1})\overrightarrow{\textup{sig}}_{\alpha_1}(v_j^*-v_i^*)}_{\text{flocking force for evaders}}\\
		\vspace{0.2cm}
		\displaystyle \qquad\quad \underbrace{-\gamma^*_i~ \overrightarrow{\textup{sig}}_{\alpha_2}(v_i-v_i^*)}_{\text{velocity-evading force}} \hspace{0.2cm}  \underbrace{- \kappa_{c,i} \phi_{c,i} (\|x_i-x_i^*\|_{\alpha_3+1}^{\alpha_3+1} )   \overrightarrow{\textup{sig}}_{\alpha_3}(x_i-x^*_i)}_{\text{position-repulsion force}},\\
		\displaystyle (x^*_i, v^*_i)\bigg|_{t=0+} =(x^*_i(0),v^*_i(0))\in \mathbb{R}^d\times \mathbb{R}^d,
	\end{cases}
\end{equation} 
where $N$, $\gamma_i$, $\gamma^*_i$, $\kappa_{r,i}$, $\kappa_{c,i}$, $\phi_{r,i}$, and $\phi_{c,i}$ are the same as for $\eqref{wolf}-\eqref{rabbit}$, $\kappa_f$ denotes the strength of relative flocking force between pursuers and evaders, disrupting force for the flocking of evaders, and flocking force for evaders, respectively, and $\phi_f$ denotes a kernel functional regarding relative flocking force between pursuers and evaders, disrupting force for the flocking of evaders, and flocking force for evaders, respectively. 
{{Motivated by the above example, the parameters $\alpha_1$, $\alpha_2$, and $\alpha_3$ in $(0,1)$ are fractional exponents that determine the non-Lipschitz structure of the relative flocking force, the velocity-tracking and velocity-evading forces, and the position-attraction and position-repulsion forces, respectively, thereby enabling finite-time convergence.}}

The term labeled as the \emph{{{disrupting}} force for the flocking of evaders} in \eqref{Fwolf} and the term labeled as the \emph{flocking force for evaders} in \eqref{Frabbit} have the same mathematical form. This term is intentionally added to the pursuer dynamics so that, later in the analysis, by subtracting the two systems, we can derive a closed system for the relative variables
\[
\hat{x}_i := x_i - x_i^*, \quad \hat{v}_i := v_i - v_i^*.
\]
Since the mapping $\overrightarrow{\mathrm{sig}}_{\alpha}$ is nonlinear, simply subtracting the pursuer and evader dynamics does not give a convenient form for the evolution of the relative variables. For this reason, we include in the pursuer dynamics an auxiliary matching term having the same form as the evader flocking interaction. In this way, after subtracting \eqref{Frabbit} from \eqref{Fwolf}, the resulting system can be written in a closed form with respect to \((\hat{x}_i,\hat{v}_i)\). Therefore, this term should be understood not as an additional independent physical interaction, but as a structural device for formulating and analyzing the relative tracking dynamics. Moreover, with this modeling choice, when \(\alpha_1=\alpha_2=\alpha_3=1\), the system \eqref{Fwolf}--\eqref{Frabbit} reduces exactly to the original asymptotic model \eqref{wolf}--\eqref{rabbit}. Thus, the present setting is naturally consistent with the original model.

Analyzing nonsmooth coupled systems such as the coupled system $\eqref{Fwolf}-\eqref{Frabbit}$ presents several challenges. First, since the systems are coupled, two distinct systems must be analyzed simultaneously, introducing significant complexity in understanding their interactions and overall behavior. Second, when dealing with nonsmooth vector fields, traditional smooth methods cannot be directly applied, necessitating the use of a nonsmooth framework, which complicates the analysis further. Moreover, the dynamic nature of these interactions and the nonsmoothness of vector fields pose difficulties in conducting finite-time or asymptotic stability analysis, making it challenging to derive analytical results.

{{We then}} offer basic concepts for the finite-time tracking success or failure of $\eqref{Fwolf}-\eqref{Frabbit}$.
\begin{definition}\label{D1.1}~\emph{(Finite-time tracking success or failure)}
	{{Let
	$
	Z^0:=\{(x_i(0),v_i(0),x_i^*(0),v_i^*(0))\}_{i=1}^N
	$
	be a given initial datum, and let $\mathscr{S}(Z^0)$ denote the set of all admissible solutions to the coupled system $\eqref{Fwolf}-\eqref{Frabbit}$ issued from $Z^0$, where the precise solution concept will be specified in Section~\ref{sec:2.1}.}}
	Let $\{ (x_i, v_i) \}_{i=1}^{N}$ and $\{ (x_i^*,v_i^*)\}_{i=1}^{N}$ be the configurations of pursuers and evaders, respectively. 
	\begin{enumerate}
		\item \emph{(Asymptotic tracking):} For $i\in [N]$, the $i$-th pursuer successfully asymptotically tracks the $i$-th evader if{{, for some $p\geq1$, every solution $Z(\cdot)\in \mathscr{S}(Z^0)$ satisfies}}
		\[
		\lim_{t \to \infty} \|x_i(t)- x_i^*(t)\|_p=0.
		\]
		\item \emph{(Finite-time tracking):} For $i\in [N]$, the $i$-th pursuer successfully tracks the $i$-th evader within a finite time if there exists a finite time {{$T=T(Z^0)\in [0,\infty)$ such that, for some $p\geq1$, every solution $Z(\cdot)\in \mathscr{S}(Z^0)$ satisfies}}
		\[
		\sup_{t \geq T} \|x_i(t)- x_i^*(t)\|_p=0.
		\]
		\item \emph{(Collision avoidance):} For $i\in [N]$, the $i$-th pursuer fails to track the $i$-th evader if {{there exists a positive constant $L_i=L_i(Z^0)>0$ such that, for some $p\geq1$, every solution $Z(\cdot)\in \mathscr{S}(Z^0)$ satisfies}}
		\[
		\inf_{t\geq 0}\|x_i(t)-x_i^*(t)\|_p{{\geq L_i>0}}.
		\]
	\end{enumerate}
\end{definition}

{
Now, we summarize the main conclusions of this paper through the following steps. These conclusions are obtained under appropriate sufficient conditions on the initial data, system parameters, and kernel functionals, which will be specified later.
	
\begin{itemize}
	\item (Step A):
	For a fixed admissible initial datum, we introduce the relative
	position and velocity variables
	\[
	\hat{x}_i:=x_i-x_i^*,
	\quad
	\hat{v}_i:=v_i-v_i^*,
	\]
	and reformulate the proposed pursuer--evader system in terms of
	$\{(\hat{x}_i,\hat{v}_i)\}_{i=1}^N$.
	
	\vspace{0.1cm}
	
	\item (Step B):
	We construct an appropriate energy functional $\mathcal{E}$ and
	derive its monotonicity along every global classical solution:
	\[
	\frac{d\mathcal{E}}{dt}\leq 0.
	\]
	This yields uniform bounds for the relative positions and velocities.
	
	\vspace{0.1cm}
	
	\item (Step C):
	Using these energy estimates together with an invariance argument,
	we prove asymptotic tracking success for the pursuers in
	$\mathcal{N}_1$, while showing that the pursuers in $\mathcal{N}_2$
	remain uniformly separated from their corresponding evaders. More
	precisely, for every global classical solution associated with the
	given initial datum,
	\[
	\lim_{t\to\infty}\hat{x}_i(t)
	=
	\lim_{t\to\infty}\hat{v}_i(t)
	=0,
	\quad i\in\mathcal{N}_1,
	\]
	and there exists a positive constant $L_i$ such that
	\[
	\inf_{t\geq0}\|\hat{x}_i(t)\|_{\alpha_3+1}
	\geq L_i>0,
	\quad i\in\mathcal{N}_2.
	\]
	
	\vspace{0.1cm}
	
	\item (Step D):
	Finally, using the results obtained in the preceding steps together
	with finite-time stability tools, we prove that for each admissible
	initial datum there exists a finite time $T=T(Z^0)\geq0$ such that
	every global classical solution issued from that initial datum
	satisfies
	\[
	\hat{x}_i(t)=\hat{v}_i(t)=0,
	\quad
	i\in\mathcal{N}_1,\quad t\geq T(Z^0).
	\]
	Thus, the pursuers in $\mathcal{N}_1$ successfully track their
	corresponding evaders within a finite time.
\end{itemize}
}

\vspace{0.12cm}

The rest of this paper is organized as follows. In Section \ref{sec:2}, we introduce several useful frameworks related to finite-time stability and the nonsmooth LaSalle invariance principle, together with some associated concepts, {{and also review}} previous results for the smooth coupled system $\eqref{wolf}-\eqref{rabbit}$. In Section \ref{sec:3}, we establish a sufficient framework for determining the success or failure of finite-time tracking in the nonsmooth coupled system $\eqref{Fwolf}-\eqref{Frabbit}$. Finally, Section \ref{sec:4} summarizes the main results and discusses open problems for future research.

	\section{Preliminaries} \label{sec:2}
	\setcounter{equation}{0}
In this section, we first present a useful lemma on finite-time stability, along with the fundamental concepts required to understand it. We then introduce the nonsmooth LaSalle invariance principle and review existing results for the coupled system $\eqref{wolf}-\eqref{rabbit}$.
	
\subsection{Finite-time stability lemma}\label{sec:2.1}
In this subsection, we introduce a key lemma concerning finite-time stability. 
First, we define the state vector $x \in \mathbb{R}^m$ and the nonsmooth vector field $F: \mathbb{R}^m \to \mathbb{R}^m$ as follows:
\begin{align*}
	x = (x_1, \ldots, x_m), \quad F = (F_1, \ldots, F_m),
\end{align*}
where $F$ denotes a nonsmooth vector field that is continuous but not necessarily locally Lipschitz continuous.

We then consider the Cauchy problem for the following first-order ODE system:
\begin{align}
	\begin{cases}\label{B0}
		\dot{x}(t) = F(x(t)), \quad t > 0, \\
		x(0) = x^0.
	\end{cases}
\end{align}
{{Since the vector field $F$ may be nonsmooth, the standard
Cauchy--Lipschitz (Picard--Lindel\"of) theorem cannot be used to obtain
local uniqueness.  The continuity of $F$ nevertheless allows us to apply
Peano's theorem, which guarantees at least one local classical solution.
More precisely, every local solution admits a maximal continuation
\[
x\in C^1([0,T_{\max});\mathbb{R}^m),
\qquad 0<T_{\max}\leq\infty.
\]
Peano's theorem alone does not imply that $T_{\max}=\infty$; global
existence must be verified separately for the system under consideration.

A solution is called \emph{global} if $T_{\max}=\infty$.  Whenever global
existence has been established for every maximal solution issued from
$x^0$, we use the notation
\[
\mathscr{S}(x^0)
:=
\left\{
 x(\cdot):x(0)=x^0,\;
 x\in C^1([0,\infty);\mathbb{R}^m),\;
 \dot x=F(x)
\right\}.
\]
Uniqueness is not assumed, and a statement concerning
$\mathscr{S}(x^0)$ is always understood to hold for every solution in
this set.  For the coupled pursuer--evader system studied below, the
required global-existence property is proved in
Corollary~\ref{C3.1}.}}

{{
We now introduce the two localized notions used below: forward positive
invariance and homogeneity on a bounded set under contractive dilations.}}

{{
\begin{definition}\label{D2.1}
\emph{\cite{B-B,L}}~\emph{(Localized notions for finite-time stability)}
Let $K\subset\mathbb{R}^m$ be a set containing the origin, and let
$(d_1,\ldots,d_m)\in(0,\infty)^m$. For $\varepsilon\in(0,1]$, define
\[
\delta_\varepsilon(x_1,\ldots,x_m)
:=
(\varepsilon^{d_1}x_1,\ldots,\varepsilon^{d_m}x_m).
\]
\begin{enumerate}
\item
The set $K$ is said to be forward positively invariant for~\eqref{B0} if
every maximal classical solution
\[
x\in C^1([0,T_{\max});\mathbb R^m),
\qquad x(0)\in K,
\]
is global and satisfies
\[
T_{\max}=\infty,
\qquad
x(t)\in K,\qquad t\geq0.
\]
Thus, forward positive invariance is understood here in the strong sense,
for every possibly non-unique maximal solution issued from $K$.

\item
The vector field $F$ is said to be homogeneous of degree
$\lambda\in\mathbb R$ on $K$ with respect to the dilation
$(d_1,\ldots,d_m)$ if, for every $x\in K$,
$\varepsilon\in(0,1]$, and $i\in\{1,\ldots,m\}$,
\[
F_i(\delta_\varepsilon x)
=
\varepsilon^{\lambda+d_i}F_i(x).
\]
\end{enumerate}
\end{definition}

With the localized notions introduced above, we now present a finite-time
stability lemma tailored to the setting of a bounded forward positively
invariant set. Before stating the lemma, we recall that the origin is said
to be asymptotically stable relative to $K$ if it is Lyapunov stable and
attractive for all global solutions of~\eqref{B0} whose initial data
belong to $K$. The following lemma may be viewed as a localized counterpart
of the classical finite-time stability results in~\cite{B-B,L}.

\begin{lemma}\label{L2.1}~\emph{(Finite-time stability under negative homogeneity)}
Let $F:\mathbb{R}^m\to\mathbb{R}^m$ be a continuous vector field
satisfying $F(0)=0$, and let $K\subset\mathbb{R}^m$ be a bounded and
forward positively invariant set containing the origin in its interior.
Suppose that $F$ is homogeneous of degree $\lambda<0$ on $K$ with
respect to a dilation $(d_1,\ldots,d_m)$. Assume that the origin is
asymptotically stable relative to $K$. Then there exist a nonempty open
neighborhood $V\subset K$ of the origin and a finite constant $T_*>0$
such that, for every initial condition $x^0\in V$, every global solution
of~\eqref{B0} satisfies
\[
x(t)=0,
\qquad t\geq T_*.
\]
Consequently, the origin is locally uniformly finite-time stable relative
to $K$.
\end{lemma}
}}

{{
\begin{proof}
For the dilation in Definition~\ref{D2.1}, introduce the homogeneous
gauge
\[
|x|_d:=\max_{1\leq i\leq m}|x_i|^{1/d_i},
\qquad
B_r^d:=\{x\in\mathbb R^m:|x|_d<r\}.
\]
Then
\[
|\delta_\varepsilon x|_d=\varepsilon|x|_d,
\qquad
\delta_\varepsilon B_r^d=B_{\varepsilon r}^d.
\]
Since the origin belongs to the interior of $K$, we may choose $R>0$
such that
\[
\overline{B_R^d}\subset\operatorname{int}K.
\]
By Lyapunov stability, there exists $r\in(0,R)$ such that every global
solution satisfying $x(0)\in\overline{B_r^d}$ remains in $B_R^d$ for all
$t\geq0$. Fix any $\rho\in(0,1)$. Applying Lyapunov stability once more,
choose an open neighborhood $W$ of the origin such that
\[
W\subset B_{\rho r}^d
\]
and every solution starting in $W$ remains in $B_{\rho r}^d$ for all
future times.

We first establish a uniform entrance time. We claim that there exists
$T>0$ such that every global solution satisfying
$x(0)\in\overline{B_r^d}$ satisfies
\[
x(t)\in B_{\rho r}^d,
\qquad t\geq T.
\]
Suppose otherwise. Then there exist global solutions $x_n(\cdot)$ with
$x_n(0)\in\overline{B_r^d}$ and times $t_n\geq n$ such that
$x_n(t_n)\notin B_{\rho r}^d$. All these solutions remain in the common
compact set $\overline{B_R^d}$. Since $F$ is bounded there, the family
$\{x_n\}$ is equicontinuous on every bounded time interval. By
Arzel\`a--Ascoli and a diagonal argument, a subsequence converges
uniformly on compact time intervals to a global solution $x_\infty$
with $x_\infty(0)\in\overline{B_r^d}$. Passing to the limit in the
integral equation verifies that $\dot x_\infty=F(x_\infty)$.
Attractivity gives $x_\infty(T_0)\in W$ for some $T_0>0$. Hence
$x_n(T_0)\in W$ for all sufficiently large $n$. The time-shifted
solutions then remain in $B_{\rho r}^d$ for all later times, contradicting
$x_n(t_n)\notin B_{\rho r}^d$. This proves the claim.

We next use negative homogeneity. Let $k\geq0$ and let $y(\cdot)$ be any
solution satisfying
\[
y(0)\in\overline{B_{\rho^k r}^d}.
\]
On the interval on which the rescaled curve is defined, set
\[
x(s):=\delta_{\rho^{-k}}y(\rho^{-k\lambda}s).
\]
Then $x(0)\in\overline{B_r^d}$. As long as $x(s)\in B_R^d\subset K$,
the homogeneity identity yields, componentwise,
\[
\frac{dx_i}{ds}
=
\rho^{-k\lambda}\rho^{-kd_i}
F_i(\delta_{\rho^k}x)
=F_i(x).
\]
Since $y$ is global, the rescaled curve $x$ is defined for all $s\geq0$.
We claim that it cannot leave $B_R^d$. Otherwise, let $s_0$ be its first
exit time. Because $\overline{B_R^d}\Subset K$, the homogeneity identity
continues to hold on a short interval beyond $s_0$; hence the restriction
of $x$ to that interval is a classical solution of~\eqref{B0}. Extend
this solution segment to a maximal classical solution of~\eqref{B0}.
Its initial value belongs to $K$, so forward positive invariance makes the
extension global, and Lyapunov stability forces it to remain in $B_R^d$.
This contradicts the definition of $s_0$. Hence the rescaled curve is a
global solution of~\eqref{B0} contained in $B_R^d$. By the uniform
entrance estimate,
$x(T)\in B_{\rho r}^d$, and therefore
\[
y(\rho^{-k\lambda}T)
=
\delta_{\rho^k}x(T)
\in B_{\rho^{k+1}r}^d.
\]
Applying this argument successively to time shifts of a solution starting
in $\overline{B_r^d}$, we find that it enters
$B_{\rho^n r}^d$ no later than
\[
t_n:=\sum_{k=0}^{n-1}\rho^{-k\lambda}T.
\]
Since $\lambda<0$ and $\rho\in(0,1)$,
\[
T_*:=\lim_{n\to\infty}t_n
=
\frac{T}{1-\rho^{-\lambda}}
<\infty.
\]
At the times $t_n$, the solution lies in $B_{\rho^n r}^d$, and hence
converges to the origin. Continuity gives $x(T_*)=0$. Finally, Lyapunov
stability with respect to all solutions implies that every solution
starting at the origin remains identically equal to zero. Thus every
solution starting in
\[
V:=B_r^d
\]
satisfies $x(t)=0$ for all $t\geq T_*$. The constants $V$ and $T_*$ are
independent of the choice of solution.
\end{proof}

}}

\subsection{Nonsmooth LaSalle's invariance principle} \label{sec:2.2}
In this subsection, we present the nonsmooth LaSalle invariance principle \cite{R-H-L}, which will be applied in Sections \ref{sec:3} and \ref{sec:4}. The classical LaSalle invariance principle, as described in \cite{L-R}, addresses the asymptotic stability of equilibria for systems with a locally Lipschitz continuous vector field. However, when considering the asymptotic or finite-time stability of equilibria for nonsmooth dynamical systems, the standard approach used in \cite{L-R} does not apply directly. To overcome this challenge, we introduce a more general version of LaSalle's invariance principle, specifically designed for nonsmooth systems.

We now reconsider the following first-order ODE system, where the vector field $F$ on $\mathbb{R}^m$ is continuous but not locally Lipschitz continuous:

\begin{align}\label{B1}
\begin{cases} 
	\dot{x} = F(x),  \quad t > 0, \\ 
	x(t)\Big|_{t = 0+} = x^0.
\end{cases}
\end{align}

{{In what follows, the positive limit set is attached to a chosen
solution, since different solutions issued from the same initial datum
may have different limit sets.}}
\begin{definition}\label{D2.2}\emph{\cite{L-R,R-H-L}}~\emph{(Positive limit set)}
	{{Let $x(\cdot)$ be a global solution of~\eqref{B1}.  Its
	positive limit set, denoted by $\omega(x(\cdot))$, is defined by}}
	\[
	{{\omega(x(\cdot)) := \left\{ y \in \mathbb{R}^m \, : \, \exists \, \{ t_n \}_{n \geq 1} \quad \text{such that} \quad \lim_{n \to \infty} t_n = \infty \quad \text{and} \quad \lim_{n \to \infty} x(t_n) = y \right\}.}}
	\]
\end{definition}

	It is well-established that for an ODE with a locally Lipschitz continuous vector field, the positive limit set is closed, positively invariant, and connected (see \cite{L-R}). However, for an ODE with a non-locally Lipschitz continuous vector field, the positive limit set is only semi-invariant. 

Next, we introduce the concept of semi-invariance for the positive limit set. 
{
\begin{definition}\label{D2.3}\emph{\cite{R-H-L}}~\emph{(Semi-invariance)}
A set $S\subset\Omega$ is said to be \emph{semi-invariant} for
\begin{align}\label{NNN-1}
\dot{x}=F(x)
\end{align}
if, for every $y\in S$, there exists a maximal (non-continuable) solution
$x:J\to\Omega$, defined on an interval $J$ containing $0$, such that
\[
x(0)=y,
\qquad
x(t)\in S\quad\text{for all }t\in J.
\]
Thus semi-invariance requires the existence of at least one trajectory
through each point of $S$ that remains in $S$; it does not require every
trajectory through that point to remain in $S$.
\end{definition}
}

	It is important to note that every positive limit set of an autonomous system is semi-invariant. Specifically, for the system \eqref{NNN-1}, the intersection $\omega(x(\cdot)) \cap \Omega$ is semi-invariant. For details, we refer to the aforementioned work.

Finally, we state the nonsmooth version of LaSalle's invariance principle as follows:
\begin{proposition}\label{P2.2}\emph{\cite{R-H-L}}~\emph{(Nonsmooth LaSalle's Invariance Principle)}  
	Let $F: \mathbb{R}^m \to \mathbb{R}^m$ be a continuous (nonsmooth) vector field and $\Omega \subset \mathbb{R}^m$ an open subset. Suppose that $x = x(t)$ is a solution to the ODE system \eqref{B1}, and let $\mathcal{L}$ be a locally Lipschitz continuous functional such that
	\[
	D^+ \mathcal{L}(x(t)) \leq 0, \quad \text{where} \quad D^+ \text{ denotes the upper Dini derivative}.
	\]
	Then, the positive limit set {{$\omega(x(\cdot)) \cap \Omega$}} is contained within the zero set of the upper Dini derivative of $\mathcal{L}$, i.e.,
	\[
	{{\omega(x(\cdot)) \cap \Omega}} \subset \{ x \in \Omega \, : \, D^+ \mathcal{L} = 0 \}.
	\]
\end{proposition}

\subsection{Previous results} \label{sec:2.3}

{
	In this subsection, we briefly summarize the main results studied in \cite{A-B-H-Y1} for the smooth coupled pursuer--evader system $\eqref{wolf}-\eqref{rabbit}$. For this purpose, the entire ensemble $[N]$ is partitioned into two disjoint, nonempty subensembles:
	\[
	\mathcal{N}_1:=\{i_1,\ldots,i_{N_1}\},
	\quad
	\mathcal{N}_2:=\{j_1,\ldots,j_{N_2}\},
	\quad
	\mathcal{N}_1\cup\mathcal{N}_2=[N],
	\quad
	\mathcal{N}_1\cap\mathcal{N}_2=\emptyset,
	\quad
	N_1,N_2\in\mathbb{N}.
	\]
	
	Under suitable conditions on the initial data, system parameters, and kernel functionals, it was shown in \cite{A-B-H-Y1} that the pursuers in $\mathcal{N}_1$ successfully track their corresponding evaders asymptotically, whereas those in $\mathcal{N}_2$ fail to do so. More precisely, one has
	\[
	\lim_{t\to\infty}\|x_i(t)-x_i^*(t)\|=0,
	\quad i\in\mathcal{N}_1,
	\]
	while there exist positive constants $L_i$ and a nonnegative constant $\underline{r}$ such that
	\[
	\inf_{t\geq0}\|x_i(t)-x_i^*(t)\|
	\geq L_i>\sqrt{\underline{r}}\geq0,
	\quad i\in\mathcal{N}_2.
	\]
	These qualitative results were obtained by using the standard LaSalle invariance principle \cite{L-R} together with appropriate energy estimates.
	
Furthermore, by combining the hypocoercivity method (see \cite{Vi}) with suitable technical estimates, the pursuers in $\mathcal{N}_1$ were shown to achieve asymptotic tracking at an exponential rate. Namely, there exist positive constants $\lambda_1$ and $\lambda_2$ such that
	\[
	\sum_{i\in\mathcal{N}_1}
	\left(
	\|x_i(t)-x_i^*(t)\|^2
	+
	\|v_i(t)-v_i^*(t)\|^2
	\right)
	\leq
	\lambda_1e^{-\lambda_2t},
	\quad t\geq0.
	\]
	For the precise sufficient conditions and detailed proofs of these results, we refer the reader to \cite{A-B-H-Y1}.
	
	The present study extends this previous smooth framework to a nonsmooth setting and develops sufficient conditions under which the pursuers in $\mathcal{N}_1$ achieve tracking within a finite time, while those in $\mathcal{N}_2$ continue to exhibit tracking failure.
}

\section{Finite-time tracking success or failure} \label{sec:3}
\setcounter{equation}{0}
In this section, we analyze the nonsmooth coupled pursuer--evader system $\eqref{Fwolf}-\eqref{Frabbit}$ to examine the success or failure of finite-time tracking. The main objective is to establish a rigorous framework for determining sufficient conditions under which each pursuer successfully tracks its evader within a finite time or fails to do so. Before delving into the main discussion, we first reformulate $\eqref{Fwolf}-\eqref{Frabbit}$ to facilitate obtaining the desired results in Section \ref{sec:3}. Specifically, for $i \in [N]$, we define  
\[\hat{x}_i := x_i - x_i^*, \quad \hat{v}_i := v_i - v_i^*,\]  
and restate $\eqref{Fwolf}-\eqref{Frabbit}$ as follows:
\begin{equation}\label{Fwolf'}
	\begin{cases} 
		\vspace{0.2cm}
		\displaystyle {\dot{x}_i} =v_i,\quad t>0,\quad  i\in [N]:=\{1,\ldots,N\},\\
		\vspace{0.2cm}
		\displaystyle \dot{v}_i={\frac{\kappa_{f}}{N}\sum_{j=1}^{N}\phi_f(\|(x_i-x_j)-(x_i^*-x_j^*)\|_{\alpha_3+1})\overrightarrow{\textup{sig}}_{\alpha_1}((v_j-v_j^*)- (v_i-v_i^*))}  \\
		\vspace{0.2cm}
		\displaystyle\qquad\quad+{\frac{\kappa_{f}}{N}\sum_{j=1}^{N}\phi_f(\|(x_i-x_j)-(x_i^*-x_j^*)\|_{\alpha_3+1})\overrightarrow{\textup{sig}}_{\alpha_1}(v_j^*-v_i^*)} \\
		\vspace{0.2cm}
		\displaystyle \qquad\quad {-\gamma_i~\overrightarrow{\textup{sig}}_{\alpha_2} (v_i-v_i^*)} {- \kappa_{r,i} \phi_{r,i} (\|x_i-x_i^*\|_{\alpha_3+1}^{\alpha_3+1} )   \overrightarrow{\textup{sig}}_{\alpha_3}(x_i-x^*_i)},\\
		\displaystyle (x_i, v_i)\bigg|_{t=0+} =(x_i(0),v_i(0))\in \mathbb{R}^d\times \mathbb{R}^d,
	\end{cases}
\end{equation}
\begin{equation}\label{Frabbit'}
	\begin{cases} 
		\vspace{0.2cm}
		\displaystyle {\dot{x}^*_i} =v^*_i,\quad t>0,\quad  i\in [N]:=\{1,\ldots,N\},\\
		\vspace{0.2cm}
		\displaystyle \dot{v}^*_i={\frac{\kappa_{f}}{N}\sum_{j=1}^{N}\phi_f(\|(x_i-x_j)-(x_i^*-x_j^*)\|_{\alpha_3+1})\overrightarrow{\textup{sig}}_{\alpha_1}(v_j^*-v_i^*)}\\
		\vspace{0.2cm}
		\displaystyle \qquad\quad {-\gamma^*_i~ \overrightarrow{\textup{sig}}_{\alpha_2}(v_i-v_i^*)} {- \kappa_{c,i} \phi_{c,i} (\|x_i-x_i^*\|_{\alpha_3+1}^{\alpha_3+1} )   \overrightarrow{\textup{sig}}_{\alpha_3}(x_i-x^*_i)},\\
		\displaystyle (x^*_i, v^*_i)\bigg|_{t=0+} =(x^*_i(0),v^*_i(0))\in \mathbb{R}^d\times \mathbb{R}^d.
	\end{cases}
\end{equation} 
Here, by subtracting \eqref{Frabbit'} from \eqref{Fwolf'}, we derive the following system for $\{(\hat{x}_i,\hat{v}_i)\}_{i=1}^N$:
\begin{equation}\label{Fw-r}
	\begin{cases} 
		\vspace{0.2cm}
		\displaystyle {\dot{\hat{x}}_i} =\hat{v}_i,\quad t>0,\quad  i\in [N],\\
		\vspace{0.2cm}
		\displaystyle \dot{\hat{v}}_i=\frac{\kappa_{f}}{N}\sum_{j=1}^{N}\phi_f(\|\hat{x}_i-\hat{x}_j\|_{\alpha_3+1})\overrightarrow{\textup{sig}}_{\alpha_1}\left(\hat{v}_j- \hat{v}_i\right) -(\gamma_i-\gamma_i^*) \overrightarrow{\textup{sig}}_{\alpha_2}(\hat{v}_i) \\
			\vspace{0.2cm}
		\displaystyle \hspace{1cm}- \kappa_{r,i} \phi_{r,i} (\|\hat{x}_i\|_{\alpha_3+1}^{\alpha_3+1} )\overrightarrow{\textup{sig}}_{\alpha_3}(\hat{x}_i)+\kappa_{c,i} \phi_{c,i} (\|\hat{x}_i\|_{\alpha_3+1}^{\alpha_3+1} )\overrightarrow{\textup{sig}}_{\alpha_3}(\hat{x}_i),\\
		\displaystyle (x_i, v_i,x^*_i, v^*_i)\bigg|_{t=0+} =(x_i(0),v_i(0),x^*_i(0),v^*_i(0))\in \mathbb{R}^{4d}.
	\end{cases}
\end{equation}
Throughout this section, to investigate the complex coupled system $\eqref{Fwolf'}-\eqref{Frabbit'}$, we employ the simplified system \eqref{Fw-r} concerning position--velocity error, $\{(\hat{x}_i,\hat{v}_i)\}_{i=1}^N$, between pursuers and evaders.

\subsection{Sufficient framework} \label{sec:3.1} 
In this subsection, we describe several sufficient frameworks $(\mathcal{G})$ for the desired finite-time tracking success or failure of $\eqref{Fwolf}-\eqref{Frabbit}$. To do this, we first divide $[N]$ into two disjoint index sets as follows: 
	\[
	\mathcal{N}_1:=\{i_1,\ldots,i_{N_1}\},
	\quad
	\mathcal{N}_2:=\{j_1,\ldots,j_{N_2}\},
	\quad
	\mathcal{N}_1{{\cup}}\mathcal{N}_2=[N],
	\quad
	{{\mathcal{N}_1\cap\mathcal{N}_2=\emptyset}},
	\quad
	N_1,N_2\in\mathbb{N},
	\]
and we develop a framework in which the pursuers in the set ${\mathcal N}_1$ successfully achieve finite-time tracking, while those in ${\mathcal N}_2$ do not.
\vspace{0.2cm}
\begin{itemize}
	\item $(\mathcal{G}1)$ (Nonnegativity of system parameters): We set
\[ \kappa_{f} \geq 0, \quad \kappa_{c,i} > 0, \quad  \kappa_{r,i} > 0, \quad \gamma_i\geq 0,\quad \gamma_i^*\geq 0,\hspace{0.2cm} \forall~i\in [N],\quad \alpha_i\in (0,1),\hspace{0.2cm} i\in [3]. \] 
	\item
	$(\mathcal{G}2)$ (Conditions for asymptotic tracking success): We assume that $\phi_f:[0,\infty) \to [0,\infty)$ is continuous functional, and for each $i\in \mathcal{N}_1$, $\phi_{r,i}:[0,\infty) \to [0,\infty)$ and $\phi_{c,i}:[0,\infty) \to [0,\infty)$ are continuous functionals. Moreover, we define
	\[\psi_{r,i}(u):=\kappa_{r,i} \phi_{r,i} (u)-\kappa_{c,i} \phi_{c,i} (u),\quad u\geq 0,\hspace{0.2cm} \forall~ i\in \mathcal{N}_1,\]
{ where each $\psi_{r,i}:[0,\infty)\to[0,\infty)$ satisfies
\[
\psi_{r,i}(0)\geq 0,
\quad
\psi_{r,i}(u)>0,\hspace{0.2cm} \forall\,u>0.
\]} Suppose that for each $i\in \mathcal{N}_1$,
		\begin{equation*}
		\begin{cases} 
			\vspace{0.2cm}
			\displaystyle ~\gamma_i>\gamma_i^*,\quad  \mathcal{E}(0) < \frac{1}{(\alpha_3+1)}\int_{0}^{\infty} \psi_{r,i} (u) du,\\
			\displaystyle~U_i:= \sup_{s\in [0,\infty)}\left\{ \int_{0}^{s^{\alpha_3+1}} \psi_{r,i}(u)du = (\alpha_3+1) \mathcal{E}(0)\right \} \in  \left[0, \infty\right),
		\end{cases}
	\end{equation*} where $\mathcal{E}(0)$ is given by
    \begin{align}\label{C0}
    	\begin{aligned}
    	\mathcal{E}(0)&:=\frac{1}{2}\sum_{i=1}^{N}\|\hat{v}_i(0)\|^2+\frac{1 }{(\alpha_3+1)}\sum_{i\in \mathcal{N}_1}\left(\int_{0}^{\|\hat{x}_i(0)\|^{\alpha_3+1}_{\alpha_3+1}} \psi_{r,i}(u) du\right)\\
    	&\qquad + \frac{1 }{(\alpha_3+1)}\sum_{i\in \mathcal{N}_2}\left(\kappa_{r,i}\int_{0}^{\|\hat{x}_i(0)\|^{\alpha_3+1}_{\alpha_3+1}} \phi_{r,i}(u) du\right)\\
    	&\qquad+\frac{1}{(\alpha_3+1)}\sum_{i\in \mathcal{N}_2}\left(\kappa_{c,i}\int_{\|\hat{x}_i(0)\|^{\alpha_3+1}_{\alpha_3+1}}^{\infty} \phi_{c,i}(u) du\right).
    	\end{aligned}
    \end{align}
	\item $(\mathcal{G}3)$ (Conditions for asymptotic tracking failure):~ 
	For each $i\in \mathcal{N}_2$ and a nonnegative constant $\underline{r}\geq 0$, let $\phi_{r,i}:[0,\infty) \to [0,\infty)$ and $\phi_{c,i}:(\underline{r}, \infty) \to [0,\infty)$ be continuous functionals satisfying
	\[
	\int_{\underline{r}+\varepsilon}^\infty \phi_{c,i}(u) du < \infty, \quad \forall~\varepsilon>0.
	\]
	We further assume that for each $i\in \mathcal{N}_2$,
	\begin{equation*}
		\begin{cases} 
			\vspace{0.1cm}
			\displaystyle  \min_{i\in\mathcal{N}_2}\|x_i(0)-x_i^*(0)\|_{\alpha_3+1}> \underline{r}^{\frac{1}{\alpha_3+1}}\geq 0,\quad \gamma_i=\gamma_i^*,\\
			\vspace{0.1cm}
			\displaystyle ~\mathcal{E}(0)  < \min \left \{ \frac{\kappa_{r,i}}{(\alpha_3+1)}\int_{0}^{\infty} \phi_{r,i} (u) du,~\frac{\kappa_{c,i}}{(\alpha_3+1)}\int_{\underline{r}}^{\infty} \phi_{c,i}(u) du \right \}, \\
			\vspace{0.1cm}
			\displaystyle~U_i:= \sup_{s\in [0,\infty)}\left\{ \int_{0}^{s^{\alpha_3+1}} \phi_{r,i}(u)du = \frac{(\alpha_3+1)\mathcal{E}(0)}{\kappa_{r,i}} \right \} \in  \left[L_i, \infty\right),\\ 
				\vspace{0.1cm}
			\displaystyle ~L_i:=\inf_{s\in [0,\infty)}\left\{\int_{s^{\alpha_3+1}}^{\infty} \phi_{c,i}(u)du = \frac{(\alpha_3+1)\mathcal{E}(0)}{\kappa_{c,i}} \right\} \in (\underline{r}^{\frac{1}{\alpha_3+1}}, \infty),
		\end{cases}
	\end{equation*}where $\mathcal{E}(0)$ is defined in \eqref{C0}.
\vspace{0.2cm}
	\item $(\mathcal{G}4)$ (Conditions for finite-time tracking success):
	Suppose that 
\[\min_{i\in \mathcal{N}_2}L_i-3\max_{i\in \mathcal{N}_1}U_i{{>}}0,\]
and $\phi_f$ has a bounded support such that 
\[\text{supp}(\phi_f) \subset\left[ 0,~ \min_{i\in \mathcal{N}_2}L_i-\max_{i\in \mathcal{N}_1}U_i\right].\]
In addition, assume that one of the following four conditions holds:

\vspace{0.2cm}

\begin{enumerate}
	\item  There exist positive constants $a>b>0$ and a negative constant $\lambda<0$ such that
	\[   b=a+\lambda,\qquad b\alpha_1=b\alpha_2=b+\lambda=a\alpha_3. \] 
	For each $i\in  \mathcal{N}_1$, suppose there exist constants $\hat{\phi}_f\geq 0$ and $\hat{\psi}_{r,i} >0$   satisfying 
	\[\hspace{2.2cm}\phi_f(u)\equiv\hat{\phi}_f\hspace{0.2cm}\mbox{for}\hspace{0.2cm} u\leq 2\max_{i\in \mathcal{N}_1}U_i,\quad \text{and} \quad \psi_{r,i}(u)\equiv\hat{\psi}_{r,i}>0\hspace{0.2cm}\mbox{for}\hspace{0.2cm} u\leq U_i^{\alpha_3+1}.\]		
	\item There exist  positive constants $a>b>0$ and $\beta>0$, and a negative constant $\lambda<0$ such that 
	\[   b=a+\lambda,\qquad b\alpha_1=b\alpha_2=b+\lambda= a\beta(\alpha_3+1)+a\alpha_3. \] 
	For each $i\in  \mathcal{N}_1$, suppose there exist a constant $\hat{\phi}_f\geq 0$  satisfying
	\[\hspace{2cm}\phi_f(u)\equiv\hat{\phi}_f\hspace{0.2cm}\mbox{for}\hspace{0.2cm} u\leq 2\max_{i\in \mathcal{N}_1}U_i,\quad \text{and} \quad \psi_{r,i}(xu)=x^\beta\psi_{r,i}(u)\hspace{0.2cm}\mbox{for}\hspace{0.2cm} u\leq U_i^{\alpha_3+1},~x\in (0,1].\]
	\item There exist  positive constants $a>b>0$ and $\sigma>0$, and a negative constant $\lambda<0$ such that 
	\[   b=a+\lambda,\qquad b\alpha_1+a\sigma =b\alpha_2=b+\lambda= a\alpha_3. \] 
For each $i\in  \mathcal{N}_1$, suppose there exist a constant $\hat{\psi}_{r,i} >0$  satisfying  
	\[\hspace{2cm}\phi_f(xu)=x^\sigma\phi_f(u)\hspace{0.2cm}\mbox{for}\hspace{0.2cm} u\leq 2\max_{i\in \mathcal{N}_1}U_i,~x\in  (0,1],\quad \text{and} \quad \psi_{r,i}(u)\equiv\hat{\psi}_{r,i}\hspace{0.2cm}\mbox{for}\hspace{0.2cm} u\leq {U}_i^{\alpha_3+1}.\]
	\item There exist positive constants $a>b>0$ and $\beta>0,~ \sigma>0$, and a negative constant $\lambda<0$ such that
	\[  b=a+\lambda,\qquad b\alpha_1+a\sigma =b\alpha_2=b+\lambda=a\beta(\alpha_3+1)+a\alpha_3. \] 
	For each $i\in  \mathcal{N}_1$, suppose that $\phi_f$ and $\psi_{r,i}$ satisfy
	\begin{align*}
		&\phi_f(xu)=x^\sigma\phi_f(u),\hspace{0.2cm} u\leq 2\max_{i\in \mathcal{N}_1}U_i,~x\in  (0,1],\\ 
		&\psi_{r,i}(xu)=x^\beta\psi_{r,i}(u),\hspace{0.2cm} u\leq {U}_i^{\alpha_3+1},~x\in  (0,1].
	\end{align*}
\end{enumerate}
\end{itemize}

\begin{remark}\label{R3.1}
	The following provides a concise explanation of the sufficient conditions summarized in $(\mathcal{G})$. 
	\begin{enumerate}
		\item
		First, $(\mathcal{G}1)$ outlines the fundamental assumptions required to determine whether finite-time tracking success or failure occurs in system~\eqref{Fwolf}--\eqref{Frabbit}. Next, $(\mathcal{G}2)$ presents the sufficient conditions for asymptotic tracking success of the pursuers in the set ${\mathcal N}_1$. These conditions indicate that asymptotic success is achieved when the velocity-tracking force dominates the velocity-evading force and the position-attraction force dominates the position-repulsion force. Subsequently, $(\mathcal{G}3)$ describes the conditions leading to finite-time tracking failure for pursuers in ${\mathcal N}_2$, showing that asymptotic tracking failure arises when the velocity-tracking force does not exceed the velocity-evading force and the position-attraction force fails to dominate the position-repulsion force.

		\item
		$(\mathcal{G}1)-(\mathcal{G}3)$ ensure that the pursuers in ${\mathcal N}_1$ asymptotically converge to their corresponding evaders, and that for each $i\in {\mathcal N}_1$, the distance measured in the $(\alpha_3+1)$-norm between the pursuer and its evader remains bounded by $U_i$. Further, for each $i\in {\mathcal N}_2$, these conditions ensure that the $(\alpha_3+1)$-norm distance between the pursuer and its evader always remains greater than or equal to $L_i$.
		
		\item
		Then, $(\mathcal{G}4)$ establishes the finite-time tracking success of the pursuers in ${\mathcal N}_1$ based on the finite-time stability framework introduced in Section~\ref{sec:2.1}. In essence, the four alternative forms of $(\mathcal{G}4)$ specify different ways of constructing a vector field whose dynamics possess the homogeneous degree required in Lemma~\ref{L2.1}.
		
		\item
		Since, for each $i\in {\mathcal N}_1$, the $(\alpha_3+1)$-norm distance between the pursuer and its evader is always bounded above by $U_i$, determining the functional form of $\phi_f$ within $2\max_{i \in \mathcal{N}_1} U_i$ as prescribed in $(\mathcal{G}4)$ ensures that $\phi_f$ acts as a constant or as a homogeneous function for all $i,j\in \mathcal{N}_1$.
		
		\item
		{{Moreover, let
		\[
		M:=\max_{i\in\mathcal{N}_1}U_i,
		\quad
		L:=\min_{i\in\mathcal{N}_2}L_i.
		\]
		As will be established later, the several energy estimates associated with $(\mathcal{G}2)$ and $(\mathcal{G}3)$ imply that, for any $i,j\in\mathcal{N}_1$,
		\[
		\|\hat{x}_i-\hat{x}_j\|_{\alpha_3+1}
		\leq
		\|\hat{x}_i\|_{\alpha_3+1}
		+
		\|\hat{x}_j\|_{\alpha_3+1}
		\leq 2M,
		\]}}
		whereas, when $i\in\mathcal{N}_1$ and $j\in\mathcal{N}_2$,
		\[
		\|\hat{x}_i-\hat{x}_j\|_{\alpha_3+1}
		\geq
		\|\hat{x}_j\|_{\alpha_3+1}
		-
		\|\hat{x}_i\|_{\alpha_3+1}
		\geq {{L-M}}.
		\]
		{{Therefore, in order for $\phi_f$ to have the prescribed form for all
		within-group distances in $\mathcal{N}_1$ and to vanish for all
		cross-group distances between $\mathcal{N}_1$ and $\mathcal{N}_2$, we require
		\[
		2M<L-M,
		\]
		which is equivalent to
		\[
		L-3M>0.
		\]
		This explains the factor $3$ appearing in $(\mathcal{G}4)$. The strict
		inequality also guarantees a nonempty transition interval between
		$2M$ and $L-M$, on which the continuous kernel $\phi_f$ can decrease
		from its prescribed within-group value to zero.}}
		Consequently, $(\mathcal{G}4)$ ensures that the flocking kernel $\phi_f$
		does not act between the two groups $\mathcal{N}_1$ and
		$\mathcal{N}_2$, but only within each of the sets $\mathcal{N}_1$ or
		$\mathcal{N}_2$. This separation condition, namely, that no flocking
		interaction occurs between $\mathcal{N}_1$ and $\mathcal{N}_2$, is
		essential for defining a closed subsystem corresponding to the agents
		in $\mathcal{N}_1$. {{Theorem~\ref{T3.2} establishes attractivity of
		the origin, while Lemma~\ref{L3.2} gives Lyapunov stability. Thus the
		origin is asymptotically stable for the closed subsystem, and
		Lemma~\ref{L2.1} yields finite-time stability.}}

	\item Finally, it is worth noting that there exist infinitely many choices of initial data, system parameters, and kernel functionals satisfying $(\mathcal{G})$, {{and a concrete example satisfying this framework is provided in the following remark.}}
\end{enumerate}
\end{remark}

{
	\begin{remark}\label{R3.2}
		\emph{(A concrete example satisfying the framework $(\mathcal{G})$)}
		We provide a concrete example showing that the framework
		$(\mathcal{G})$ is nonempty.
		Let $d=2$, $N=3$, and consider the nontrivial partition
		\[
		\mathcal{N}_1=\{1,2\},
		\quad
		\mathcal{N}_2=\{3\}.
		\]
		We choose
		\[
		\alpha_1=\alpha_2=\frac12,
		\quad
		\alpha_3=\frac13,
		\quad
		a=3,
		\quad
		b=2,
		\quad
		\lambda=-1.
		\]
		Then, we have
		\[
		b=a+\lambda,
		\quad
		b\alpha_1=b\alpha_2=b+\lambda=a\alpha_3=1,
		\]
		and hence the exponent relations in $(\mathcal{G}4)(1)$ are satisfied.
			We set all coupling strengths equal to one:
		\[
		\kappa_f=\kappa_{r,i}=\kappa_{c,i}=1,
		\quad i\in[3],
		\]
		and choose the damping parameters as
		\[
		\gamma_1=\gamma_2=1,
		\quad
		\gamma_1^*=\gamma_2^*=0,
		\quad
		\gamma_3=\gamma_3^*=1.
		\]
		Thus, $\gamma_i>\gamma_i^*$ for $i\in\mathcal{N}_1$, while
		$\gamma_3=\gamma_3^*$.
		For $i\in\mathcal{N}_1$, let
		\[
		\phi_{r,i}(u)\equiv1,
		\quad
		\phi_{c,i}(u)\equiv0,
		\quad u\geq0.
		\]
		Then, one has
		\[
		\psi_{r,i}(u)
		=
		\phi_{r,i}(u)-\phi_{c,i}(u)
		\equiv1,
		\quad i\in\mathcal{N}_1.
		\]
		For the agent in $\mathcal{N}_2$, we set $\underline r=0$ and define
		\[
		\phi_{r,3}(u)
		:=
		\begin{cases}
			\displaystyle \frac14(u-3),
			& 3\leq u\leq4,\\[0.15cm]
			\displaystyle \frac14(5-u),
			& 4<u\leq5,\\[0.15cm]
			0,
			& \text{otherwise},
		\end{cases}
		\]
		and
		\[
		\phi_{c,3}(u)
		:=
		\begin{cases}
			\displaystyle \frac14(2-u),
			& 0\leq u\leq2,\\[0.15cm]
			0,
			& u>2.
		\end{cases}
		\]
		These kernels are continuous and nonnegative, and satisfy
		\[
		\int_0^\infty\phi_{r,3}(u)\,du=\frac14,
		\quad
		\int_0^\infty\phi_{c,3}(u)\,du=\frac12.
		\]
		In what follows, we choose the initial data
		\[
		x_i^*(0)=0,
		\quad
		v_i(0)=v_i^*(0)=0,
		\quad i\in[3],
		\]
		and
		\[
		x_1(0)=\left(\frac18,0\right),
		\quad
		x_2(0)=\left(-\frac18,0\right),
		\quad
		x_3(0)=\left(\left(\frac52\right)^{\frac34},0\right).
		\]
		Since $\alpha_3+1=\frac43$, it is easy to check that
		\[
		\|x_1(0)-x_1^*(0)\|_{\frac43}^{\frac43}
		=
		\|x_2(0)-x_2^*(0)\|_{\frac43}^{\frac43}
		=
		\frac1{16},
		\]
		and
		\[
		\|x_3(0)-x_3^*(0)\|_{\frac43}^{\frac43}
		=
		\frac52.
		\]
		Moreover, since one has
		\[
		\operatorname{supp}(\phi_{c,3})\subset[0,2]
		\quad\text{and}\quad
		\operatorname{supp}(\phi_{r,3})\subset[3,5],
		\]
		the contribution of the third agent to the initial energy vanishes.
		Therefore, we get
		\[
		\mathcal{E}(0)
		=
		\frac{1}{\alpha_3+1}
		\left(\frac1{16}+\frac1{16}\right)
		=
		\frac{3}{32},
		\]
		and hence
		\[
		(\alpha_3+1)\mathcal{E}(0)=\frac18.
		\]
		For $i\in\mathcal{N}_1$, the definition of $U_i$ yields
		\[
		U_1=U_2
		=
		\left(\frac18\right)^{\frac34}
		=
		2^{-\frac94}
		=:M.
		\]
		For the third agent, we observe that
		\[
		\int_0^4\phi_{r,3}(u)\,du=\frac18,
		\]
		and
		\[
		\int_1^\infty\phi_{c,3}(u)\,du=\frac18.
		\]
		Consequently, we find
		\[
		U_3=4^{\frac34}=2\sqrt2,
		\quad
		L_3=1,
		\]
		and therefore
		\[
		U_3\geq L_3>0.
		\]
		Furthermore, one can see
		\[
		\mathcal{E}(0)=\frac{3}{32}
		<
		\min\left\{
		\frac{3}{4}\int_0^\infty\phi_{r,3}(u)\,du,
		\frac{3}{4}\int_0^\infty\phi_{c,3}(u)\,du
		\right\}
		=
		\min\left\{\frac{3}{16},\frac38\right\}.
		\]
		Thus, $(\mathcal{G}2)$ and $(\mathcal{G}3)$ are satisfied.
		Finally, let
		\[
		L:=L_3=1,
		\quad
		M:=2^{-\frac94}.
		\]
		Since the following holds:
		\[
		L-3M
		=
		1-3\cdot2^{-\frac94}
		>0,
		\]
		we define the continuous flocking kernel
		\[
		\phi_f(u)
		:=
		\begin{cases}
			1,
			& 0\leq u\leq2M,\\[0.15cm]
			\displaystyle
			\frac{L-M-u}{L-3M},
			& 2M<u<L-M,\\[0.25cm]
			0,
			& u\geq L-M.
		\end{cases}
		\]
		Then, we observe
		\[
		\operatorname{supp}(\phi_f)\subset[0,L-M],
		\]
		and
		\[
		\phi_f(u)\equiv1
		\quad\text{for }0\leq u\leq2M.
		\]
		Since we have
		\[
		\psi_{r,i}(u)\equiv1,
		\quad
		i\in\mathcal{N}_1,
		\]
	    the	condition $(\mathcal{G}4)(1)$ holds with
		\[
		\hat{\phi}_f=1,
		\quad
		\hat{\psi}_{r,1}
		=
		\hat{\psi}_{r,2}
		=
		1.
		\]
		In conclusion, the framework $(\mathcal{G})$ is nonempty.
	\end{remark}
}

\subsection{Asymptotic tracking success or tracking failure}\label{sec:3.2} 
In this subsection, we examine the asymptotic tracking success of $\{(x_i, v_i)\}_{i \in \mathcal{N}_1}$ and the tracking failure of $\{(x_i, v_i)\}_{i \in \mathcal{N}_2}$ in the nonsmooth coupled system $\eqref{Fwolf}-\eqref{Frabbit}$ under the framework $(\mathcal{G}1)-(\mathcal{G}3)$. To do this, we introduce a Lyapunov functional tailored to this coupled system. With the help of this Lyapunov functional, we derive several energy estimates, including uniform upper bounds for $\|x_i - x_i^*\|_p$ and $\|v_i - v_i^*\|_p$ for some $p\geq 1$. Now, the kinetic energy $\mathcal{E}_k$, relaxation energy $\mathcal{E}_r$, potential energy $\mathcal{E}_p$, and total energy $\mathcal{E}$ are defined as follows for $t \geq 0$:
\begin{align*}
	\begin{aligned} 
		& \mathcal{E}_k(t):=\frac{1}{2}\sum_{i=1}^{N}\|\hat{v}_i(t)\|^2, \\
		& \mathcal{E}_{r}(t) :=  \frac{1 }{(\alpha_3+1)}\sum_{i\in \mathcal{N}_1}\left(\int_{0}^{\|\hat{x}_i(t)\|^{\alpha_3+1}_{\alpha_3+1}} \psi_{r,i}(u) du\right)+ \frac{1 }{(\alpha_3+1)}\sum_{i\in \mathcal{N}_2}\left(\kappa_{r,i}\int_{0}^{\|\hat{x}_i(t)\|^{\alpha_3+1}_{\alpha_3+1}} \phi_{r,i}(u) du\right), \\
		&\mathcal{E}_ {p}(t):=\frac{1}{(\alpha_3+1)}\sum_{i\in \mathcal{N}_2}\left(\kappa_{c,i}\int_{\|\hat{x}_i(t)\|^{\alpha_3+1}_{\alpha_3+1}}^{\infty} \phi_{c,i}(u) du\right),\\
		& \mathcal{E}(t):=\mathcal{E}_k(t)+\mathcal{E}_r(t) +\mathcal{E}_p(t).
	\end{aligned}
\end{align*}

We then investigate the monotonicity of $\mathcal{E}$, assuming $(\mathcal{G}1)-(\mathcal{G}3)$.
\begin{lemma}\label{L3.1}
	Let
			{{$
			Z(\cdot)
			=
			\{(x_i(t),v_i(t),x_i^*(t),v_i^*(t))\}_{i=1}^N
			$}}
			be {{any}} {{maximal}} {{classical}} solution to
			$\eqref{Fwolf}-\eqref{Frabbit}$ {{on $[0,T_{\max})$}} satisfying
			$(\mathcal{G}1)-(\mathcal{G}3)$. {{Then, for every}}
			{{$t\in(0,T_{\max})$}},
	\begin{align*}
		\frac{d \mathcal{E}}{dt} = -\frac{\kappa_f}{2N}\sum_{i,j=1}^{N}\phi_f(\|\hat{x}_i-\hat{x}_j\|_{\alpha_3+1})\left\|\hat{v}_j- \hat{v}_i\right\|^{\alpha_1+1}_{\alpha_1+1} -\sum_{i\in \mathcal{N}_1}(\gamma_i-\gamma_i^*) \|\hat{v}_i\|^{\alpha_2+1}_{\alpha_2+1}\leq 0.
   \end{align*}
\end{lemma}
\begin{proof}
	{{Fix any {{maximal}} classical solution
	\[
	Z(\cdot)
	=
	\{(x_i(t),v_i(t),x_i^*(t),v_i^*(t))\}_{i=1}^N
	{{\quad\text{on }[0,T_{\max}).}}
	\]}}
	We now differentiate $\mathcal{E}_k$ with respect to $t$ in $\eqref{Fw-r}_2$ to find
    \begin{align*}
    	\frac{d \mathcal{E}_k}{dt}=&\frac{\kappa_{f}}{N}\sum_{i,j=1}^{N}\phi_f(\|\hat{x}_i-\hat{x}_j\|_{\alpha_3+1})\left\langle \overrightarrow{\textup{sig}}_{\alpha_1}(\hat{v}_j- \hat{v}_i),\hat{v}_i\right\rangle -\sum_{i=1}^N(\gamma_i-\gamma_i^*) \langle \overrightarrow{\textup{sig}}_{\alpha_2}(\hat{v}_i), \hat{v}_i \rangle  \\
    	&-\sum_{i=1}^N \kappa_{r,i} \phi_{r,i} (\|\hat{x}_i\|^{\alpha_3+1}_{\alpha_3+1} )\langle  \overrightarrow{\textup{sig}}_{\alpha_3}(\hat{x}_i),\hat{v}_i \rangle+\sum_{i=1}^N\kappa_{c,i} \phi_{c,i} (\|\hat{x}_i\|^{\alpha_3+1}_{\alpha_3+1} )\langle \overrightarrow{\textup{sig}}_{\alpha_3}(\hat{x}_i),\hat{v}_i \rangle\\
    	=:&~\mathcal{I}_1+\mathcal{I}_2+\mathcal{I}_3+\mathcal{I}_4.
    \end{align*}
Below, we estimate $\mathcal{I}_i$ for $i \in [4]$, one by one.

\vspace{0.1cm}

\noindent $\bullet$~(Estimate of $\mathcal{I}_1$): We employ $(\mathcal{G}1)$ and the standard trick of interchanging $i$ and $j$ and dividing by $2$ to see

 \begin{align*}
	\mathcal{I}_1&=\frac{\kappa_{f}}{N}\sum_{i,j=1}^{N}\phi_f(\|\hat{x}_i-\hat{x}_j\|_{\alpha_3+1})\left\langle \overrightarrow{\textup{sig}}_{\alpha_1}(\hat{v}_j- \hat{v}_i),\hat{v}_i\right\rangle\\
	&=\frac{\kappa_{f}}{N}\sum_{i,j=1}^{N}\phi_f(\|\hat{x}_i-\hat{x}_j\|_{\alpha_3+1})\left\langle \overrightarrow{\textup{sig}}_{\alpha_1}(\hat{v}_i- \hat{v}_j),\hat{v}_j\right\rangle\\
	&=-\frac{\kappa_{f}}{2N}\sum_{i,j=1}^{N}\phi_f(\|\hat{x}_i-\hat{x}_j\|_{\alpha_3+1})\left\langle \overrightarrow{\textup{sig}}_{\alpha_1}(\hat{v}_i- \hat{v}_j),\hat{v}_i-\hat{v}_j\right\rangle\\
	&=-\frac{\kappa_f}{2N}\sum_{i,j=1}^{N}\phi_f(\|\hat{x}_i-\hat{x}_j\|_{\alpha_3+1})\left\|\hat{v}_j- \hat{v}_i\right\|^{\alpha_1+1}_{\alpha_1+1}.
\end{align*}

\noindent $\bullet$~(Estimate of $\mathcal{I}_2$): We use $(\mathcal{G}2)$ and $(\mathcal{G}3)$ to get
 \begin{align*}
	\mathcal{I}_2&=-\sum_{i=1}^N(\gamma_i-\gamma_i^*) \langle \overrightarrow{\textup{sig}}_{\alpha_2}(\hat{v}_i), \hat{v}_i \rangle\\
	&=-\sum_{i=1}^N(\gamma_i-\gamma_i^*) \|\hat{v}_i\|_{\alpha_2+1}^{\alpha_2+1}\\
	&=-\sum_{i\in \mathcal{N}_1}(\gamma_i-\gamma_i^*) \|\hat{v}_i\|_{\alpha_2+1}^{\alpha_2+1}.
\end{align*}

\noindent $\bullet$~(Estimate of $\mathcal{I}_3+\mathcal{I}_4$): It follows by a direct calculation, $(\mathcal{G}2)$, and $(\mathcal{G}3)$ that
 \begin{align*}
	\mathcal{I}_3+\mathcal{I}_4&=-\sum_{i=1}^N \kappa_{r,i} \phi_{r,i} (\|\hat{x}_i\|_{\alpha_3+1}^{\alpha_3+1} )\langle \overrightarrow{\textup{sig}}_{\alpha_3}(\hat{x}_i),\hat{v}_i \rangle+\sum_{i=1}^N\kappa_{c,i} \phi_{c,i} (\|\hat{x}_i\|^{\alpha_3+1}_{\alpha_3+1} )\langle \overrightarrow{\textup{sig}}_{\alpha_3}(\hat{x}_i),\hat{v}_i \rangle\\
	&=-\sum_{i\in \mathcal{N}_1} \left[\kappa_{r,i} \phi_{r,i} (\|\hat{x}_i\|^{\alpha_3+1}_{\alpha_3+1} )-\kappa_{c,i} \phi_{c,i} (\|\hat{x}_i\|^{\alpha_3+1}_{\alpha_3+1} )\right]\langle \overrightarrow{\textup{sig}}_{\alpha_3}(\hat{x}_i),\hat{v}_i \rangle\\
	&\quad -\sum_{i\in \mathcal{N}_2} \kappa_{r,i} \phi_{r,i} (\|\hat{x}_i\|^{\alpha_3+1}_{\alpha_3+1} )\langle \overrightarrow{\textup{sig}}_{\alpha_3}(\hat{x}_i),\hat{v}_i \rangle+\sum_{i\in\mathcal{N}_2}\kappa_{c,i} \phi_{c,i} (\|\hat{x}_i\|^{\alpha_3+1}_{\alpha_3+1} )\langle \overrightarrow{\textup{sig}}_{\alpha_3}(\hat{x}_i),\hat{v}_i \rangle\\
	&=-\sum_{i\in \mathcal{N}_1} \psi_{r,i} (\|\hat{x}_i\|^{\alpha_3+1}_{\alpha_3+1} )\langle \overrightarrow{\textup{sig}}_{\alpha_3}(\hat{x}_i),\hat{v}_i \rangle -\sum_{i\in \mathcal{N}_2} \kappa_{r,i} \phi_{r,i} (\|\hat{x}_i\|^{\alpha_3+1}_{\alpha_3+1} )\langle \overrightarrow{\textup{sig}}_{\alpha_3}(\hat{x}_i),\hat{v}_i \rangle\\
	&\quad +\sum_{i\in\mathcal{N}_2}\kappa_{c,i} \phi_{c,i} (\|\hat{x}_i\|^{\alpha_3+1}_{\alpha_3+1} )\langle \overrightarrow{\textup{sig}}_{\alpha_3}(\hat{x}_i),\hat{v}_i \rangle.
\end{align*} 
Because the following relation holds for a differentiable function $t \mapsto y(t)$:
\begin{align*}
	\left\langle \overrightarrow{\textup{sig}}_{\alpha_3}(y),\frac{dy}{dt}\right\rangle=\frac{1}{\alpha_3+1}\frac{d}{dt}\left(\|y\|_{\alpha_3+1}^{\alpha_3+1}\right),~~t > 0,
\end{align*}
it is straightforward to check
\begin{align*}
	\mathcal{I}_3+\mathcal{I}_4&=-\frac{1}{(\alpha_3+1)}\frac{d}{dt}\left(\sum_{i\in \mathcal{N}_1}\int_{0}^{\|\hat{x}_i\|^{\alpha_3+1}_{\alpha_3+1}}\psi_{r,i}(u) du\right)-\frac{1}{(\alpha_3+1)}\frac{d}{dt}\left(\sum_{i\in \mathcal{N}_2}\kappa_{r,i}\int_{0}^{\|\hat{x}_i\|^{\alpha_3+1}_{\alpha_3+1}}\phi_{r,i}(u) du\right)\\
	&\quad -\frac{1}{(\alpha_3+1)}\frac{d}{dt}\left(\sum_{i\in \mathcal{N}_2}\kappa_{c,i}\int_{\|\hat{x}_i(t)\|^{\alpha_3+1}_{\alpha_3+1}}^{\infty} \phi_{c,i}(u) du\right)\\
	&=-\frac{d\mathcal{E}_r}{dt} -\frac{d\mathcal{E}_p}{dt}.
\end{align*} Thus, we combine the estimates from $\mathcal{I}_1$ to $\mathcal{I}_4$ to yield
	\begin{align*}
	\frac{d \mathcal{E}}{dt} = -\frac{\kappa_f}{2N}\sum_{i,j=1}^{N}\phi_f(\|\hat{x}_i-\hat{x}_j\|_{\alpha_3+1})\left\|\hat{v}_j- \hat{v}_i\right\|^{\alpha_1+1}_{\alpha_1+1} -\sum_{i\in \mathcal{N}_1}(\gamma_i-\gamma_i^*) \|\hat{v}_i\|^{\alpha_2+1}_{\alpha_2+1}\leq 0.
\end{align*}
{{Since the maximal classical solution was arbitrary, the desired lemma follows on its whole interval of existence.}}
\end{proof}

{
\begin{corollary}\label{C3.1}~\emph{(Global existence of classical solutions)}
Suppose that $(\mathcal{G}1)-(\mathcal{G}3)$ hold.  Then every maximal
classical solution to the coupled system~\eqref{Fwolf}--\eqref{Frabbit}
is global.  In particular, for every admissible initial datum $Z^0$, the
solution set $\mathscr{S}(Z^0)$ is nonempty.
\end{corollary}

\begin{proof}
Let $Z(\cdot)$ be a maximal classical solution on $[0,T_{\max})$.
The energy calculation in Lemma~\ref{L3.1} is local in time and therefore
holds throughout this maximal interval. Hence
\[
\mathcal E(t)\leq\mathcal E(0),
\qquad 0\leq t<T_{\max}.
\]
The kinetic part gives
\[
\|\hat v_i(t)\|\leq\sqrt{2\mathcal E(0)},
\qquad i\in[N].
\]
The nonnegative potential terms, together with the definitions of $U_i$
and $L_i$ in $(\mathcal G2)-(\mathcal G3)$, yield
\[
\|\hat x_i(t)\|_{\alpha_3+1}\leq U_i,
\qquad i\in[N],
\]
and
\[
\|\hat x_i(t)\|_{\alpha_3+1}\geq L_i
>\underline r^{1/(\alpha_3+1)},
\qquad i\in\mathcal N_2.
\]
Thus the relative state stays in a compact subset of the domain of the
relative vector field. In particular, the arguments of all kernels remain
in compact intervals, and the boundary
$\|\hat x_i\|_{\alpha_3+1}^{\alpha_3+1}=\underline r$ of the repulsive
kernel for $i\in\mathcal N_2$ is uniformly avoided.

It remains to control the absolute motion. Let
\[
V^*(t):=\max_{i\in[N]}\|v_i^*(t)\|.
\]
For $0<\alpha<1$, norm equivalence in $\mathbb R^d$ gives a constant
$C_d>0$ such that
\[
\|\overrightarrow{\mathrm{sig}}_{\alpha}(w)\|
\leq C_d\|w\|^{\alpha}.
\]
The relative-state bounds and the continuity of the kernels therefore
imply, for every $i\in[N]$,
\[
\|\dot v_i^*(t)\|
\leq
C_0+C_1\max_{j\in[N]}\|v_j^*(t)-v_i^*(t)\|^{\alpha_1}
\leq
C\bigl(1+(V^*(t))^{\alpha_1}\bigr),
\]
where $C$ depends only on the initial energy and the fixed parameters.
Using the standard upper-Dini-derivative estimate for the maximum of
finitely many absolutely continuous functions, we obtain
\[
D^+V^*(t)
\leq C\bigl(1+(V^*(t))^{\alpha_1}\bigr)
\leq C\bigl(1+V^*(t)\bigr),
\qquad t<T_{\max}.
\]
Gronwall's inequality shows that all evader velocities are bounded on
every finite time interval. Their positions are then bounded by
integration, and the identities
$v_i=v_i^*+\hat v_i$ and $x_i=x_i^*+\hat x_i$ give the same finite-time
bounds for the pursuers.

If $T_{\max}<\infty$, the full state remains in a compact subset of the
open domain of the continuous vector field as $t\uparrow T_{\max}$.
The continuation theorem for continuous ODEs extends the solution beyond
$T_{\max}$, contradicting maximality. Hence $T_{\max}=\infty$.
Peano's theorem supplies at least one local solution, so
$\mathscr S(Z^0)\neq\emptyset$.
\end{proof}
}

 Next, with the help of Lemma \ref{L3.1}, we provide several uniform upper bounds for $\|x_i - x_i^*\|_p$'s and $\|v_i - v_i^*\|_p$'s for some $p\geq 1$.
 \begin{proposition}\label{P3.1}~\emph{(Energy estimates)}
 	Let
 	{{$
 	Z(\cdot)
 	=
 	\{(x_i(t),v_i(t),x_i^*(t),v_i^*(t))\}_{i=1}^N
 	\in \mathscr S(Z^0)
 	$}}
 	be {{any classical}} solution to the coupled system
 	$\eqref{Fwolf}-\eqref{Frabbit}$ satisfying the framework
 	$(\mathcal{G}1)-(\mathcal{G}3)$. Then, we obtain the following uniform estimates: For each $i\in [N]$,
 	\begin{align*}
 	\sup_{t\geq 0}\|v_i(t)-v_i^*(t)\|=\sup_{t\geq 0}\|\hat{v}_i(t)\|\leq \sqrt{2\mathcal{E}(0)}	,\qquad \sup_{t\geq 0}\|x_i(t)-x_i^*(t)\|_{\alpha_3+1}=\sup_{t\geq 0}\|\hat{x}_i(t)\|_{\alpha_3+1}\leq U_i.
 	\end{align*}
 \end{proposition}

\begin{proof}
	{{Fix any classical solution
		\[
		Z(\cdot)
		=
		\{(x_i(t),v_i(t),x_i^*(t),v_i^*(t))\}_{i=1}^N
		\in \mathscr S(Z^0).
		\]}}
$\bullet$~(Proof of the first assertion): Due to Lemma \ref{L3.1}, one can show that
\[ \mathcal{E}_k(t)\leq \mathcal{E}(t)\leq \mathcal{E}(0).\]
Hence, this estimate implies
\[\sup_{t\geq 0}\|v_i(t)-v_i^*(t)\|=\sup_{t\geq 0}\|\hat{v}_i(t)\|\leq \sqrt{{2\mathcal{E}(0)}},\quad \forall ~i\in [N].\]
{{Since the choice of $Z(\cdot)\in\mathscr S(Z^0)$ is arbitrary, we verify the desired first assertion.}}

\vspace{0.2cm}

\noindent $\bullet$~(Proof of the second assertion): It follows from Lemma \ref{L3.1} that
\begin{align*}
	\mathcal{E}_{r}(t) &=  \frac{1}{(\alpha_3+1)}\sum_{i\in \mathcal{N}_1}\left(\int_{0}^{\|\hat{x}_i(t)\|^{\alpha_3+1}_{\alpha_3+1}} \psi_{r,i}(u) du\right)+ \frac{1 }{(\alpha_3+1)}\sum_{i\in \mathcal{N}_2}\left(\kappa_{r,i}\int_{0}^{\|\hat{x}_i(t)\|^{\alpha_3+1}_{\alpha_3+1}} \phi_{r,i}(u) du\right)\\
	&\leq \mathcal{E}(t)\leq \mathcal{E}(0).
\end{align*}
Therefore, we have
\begin{align*}
&\frac{1}{(\alpha_3+1)}\left(\int_{0}^{\|\hat{x}_i(t)\|^{\alpha_3+1}_{\alpha_3+1}} \psi_{r,i}(u) du\right)\leq \mathcal{E}(0),\hspace{0.2cm} \forall ~i\in\mathcal{N}_1,\\ &\frac{\kappa_{r,i}}{(\alpha_3+1)}\left(\int_{0}^{\|\hat{x}_i(t)\|^{\alpha_3+1}_{\alpha_3+1}} \phi_{r,i}(u) du\right)\leq \mathcal{E}(0), \hspace{0.2cm} \forall ~i\in\mathcal{N}_2.
\end{align*} Then, for each $i\in\mathcal{N}_1$, we use the definition of $U_i$ in $({\mathcal G}2)$ to yield
\[\sup_{t\geq 0}\|x_i(t)-x_i^*(t)\|_{\alpha_3+1}=\sup_{t\geq 0}\|\hat{x}_i(t)\|_{\alpha_3+1}\leq U_i.\]
Similarly, for each $i\in\mathcal{N}_2$, Lemma \ref{L3.1} and $(\mathcal{G}3)$ lead to
\[\sup_{t\geq 0}\|x_i(t)-x_i^*(t)\|_{\alpha_3+1}=\sup_{t\geq 0}\|\hat{x}_i(t)\|_{\alpha_3+1}\leq U_i.\]
{{Since the choice of $Z(\cdot)\in\mathscr S(Z^0)$ is arbitrary, we prove the desired second assertion.}}
\end{proof}
 
Under $(\mathcal{G}1) - (\mathcal{G}3)$, we establish that the pursuer configuration $\{(x_i, v_i)\}_{i \in \mathcal{N}_1}$ successfully achieves asymptotic tracking of the evader configuration $\{(x^*_i, v^*_i)\}_{i \in \mathcal{N}_1}$. Conversely, the pursuer configuration $\{(x_i, v_i)\}_{i \in \mathcal{N}_2}$ fails to track the evader configuration $\{(x^*_i, v^*_i)\}_{i \in \mathcal{N}_2}$. To demonstrate this, we first provide a proof for the tracking failure of $\{(x_i, v_i)\}_{i \in \mathcal{N}_2}$ as follows:

 \begin{theorem}\label{T3.1}~\emph{(Tracking failure)}
		Let
			{{$
			Z(\cdot)
			=
			\{(x_i(t),v_i(t),x_i^*(t),v_i^*(t))\}_{i=1}^N
			\in \mathscr S(Z^0)
			$}}
			be {{any classical}} solution to the coupled system
			$\eqref{Fwolf}-\eqref{Frabbit}$ satisfying the framework
			$(\mathcal{G}1)-(\mathcal{G}3)$. Then, the pursuer configuration $\{(x_i, v_i)\}_{i \in \mathcal{N}_2}$ fails to track the evader configuration $\{(x^*_i, v^*_i)\}_{i \in \mathcal{N}_2}$, i.e., the following collision avoidance estimate holds: For each $i\in \mathcal{N}_2$,
	\begin{align*}
	 \inf_{t\geq 0}\|x_i(t)-x_i^*(t)\|_{\alpha_3+1}=\inf_{t\geq 0}\|\hat{x}_i(t)\|_{\alpha_3+1}\geq L_i > \underline{r}^{\frac{1}{\alpha_3+1}}  \geq 0.
	\end{align*}
\end{theorem}

\begin{proof}
{	{Fix any global classical solution
		\[
		Z(\cdot)
		=
		\{(x_i(t),v_i(t),x_i^*(t),v_i^*(t))\}_{i=1}^N
		\in \mathscr S(Z^0).
		\]}}
We apply Lemma \ref{L3.1} to deduce
	\[\mathcal{E}_ {p}(t):=\frac{1}{(\alpha_3+1)}\sum_{i\in \mathcal{N}_2}\left(\kappa_{c,i}\int_{\|\hat{x}_i(t)\|^{\alpha_3+1}_{\alpha_3+1}}^{\infty} \phi_{c,i}(u) du\right)\leq \mathcal{E}(t)\leq \mathcal{E}(0).\]
	Thus, one has
	\begin{align*}
	\frac{\kappa_{c,i}}{(\alpha_3+1)}\left(\int_{\|\hat{x}_i(t)\|^{\alpha_3+1}_{\alpha_3+1}}^\infty \phi_{c,i}(u) du\right)\leq \mathcal{E}(0), \hspace{0.2cm} \forall ~i\in\mathcal{N}_2.
	\end{align*} Hence, for each $i\in\mathcal{N}_2$, we employ the definition of $L_i$ in $({\mathcal G}3)$ to induce
	\[\inf_{t\geq 0}\|x_i(t)-x_i^*(t)\|_{\alpha_3+1}=\inf_{t\geq 0}\|\hat{x}_i(t)\|_{\alpha_3+1}\geq L_i > \underline{r}^{\frac{1}{\alpha_3+1}} \geq 0.\]
	{{Since the choice of $Z(\cdot)\in\mathscr S(Z^0)$ is arbitrary, the
	desired theorem follows.}}
\end{proof}

Now, under $(\mathcal{G}1)-(\mathcal{G}3)$, we prove that the pursuer configuration $\{(x_i, v_i)\}_{i \in \mathcal{N}_1}$ achieves asymptotic tracking success for the evader configuration $\{(x^*_i, v^*_i)\}_{i \in \mathcal{N}_1}$.
 \begin{theorem}\label{T3.2}~\emph{(Asymptotic tracking success)}
	Let
		{{$
		Z(\cdot)
		=
		\{(x_i(t),v_i(t),x_i^*(t),v_i^*(t))\}_{i=1}^N
		\in \mathscr S(Z^0)
		$}}
		be {{any global classical}} solution to the coupled system
		$\eqref{Fwolf}-\eqref{Frabbit}$ satisfying the framework
		$(\mathcal{G}1)-(\mathcal{G}3)$. Then, the pursuer configuration $\{(x_i, v_i)\}_{i \in \mathcal{N}_1}$ successfully asymptotically tracks the evader configuration $\{(x^*_i, v^*_i)\}_{i \in \mathcal{N}_1}$, i.e.,
	\begin{align*}
	\lim_{t \to \infty} (x_i(t)- x_i^*(t))=\lim_{t \to \infty} \hat{x}_i(t) =0,\quad \forall ~i\in \mathcal{N}_1.
	\end{align*} Moreover, one can derive
    	\begin{align*}
    	\lim_{t \to \infty} (v_i(t)- v_i^*(t))=\lim_{t \to \infty}\hat{v}_i(t)=0,\quad \forall ~i\in \mathcal{N}_1.
    \end{align*} 
\end{theorem}

\begin{proof}
{
	\noindent $\bullet$~(Proof of $\lim_{t \to \infty}\hat{v}_i(t)=0$):
	Let
	\[
	\hat{z}_i(t):=(\hat{x}_i(t),\hat{v}_i(t))\in \mathbb{R}^d\times \mathbb{R}^d,
	\quad i\in [N],
	\]
	and define
	\[
	\hat{Z}(t):=(\hat{z}_1(t),\ldots,\hat{z}_N(t)),
	\quad
	\hat{Z}^0:=\hat{Z}(0).
	\]
	{{Let $\omega(\hat{Z}(\cdot))$ denote the positive limit set
	(see Definition~\ref{D2.2}) of the chosen relative trajectory
	$\hat Z(\cdot)$.}}
	Proposition~\ref{P3.1} shows that the relative trajectory is precompact,
	so its positive limit set is nonempty. By Lemma~\ref{L3.1},
	\[
	\frac{d\mathcal E}{dt}
	=
	-\frac{\kappa_f}{2N}\sum_{i,j=1}^N
	\phi_f(\|\hat x_i-\hat x_j\|_{\alpha_3+1})
	\|\hat v_j-\hat v_i\|_{\alpha_1+1}^{\alpha_1+1}
	-\sum_{i\in\mathcal N_1}(\gamma_i-\gamma_i^*)
	\|\hat v_i\|_{\alpha_2+1}^{\alpha_2+1}.
	\]
	Since $\gamma_i-\gamma_i^*>0$ for $i\in\mathcal N_1$, the nonsmooth
	LaSalle principle implies that every point of
	$\omega(\hat Z(\cdot))$ satisfies
	\[
	\hat v_i=0,
	\quad i\in\mathcal N_1,
	\]
	and
	\[
	\kappa_f\phi_f(\|\hat x_i-\hat x_j\|_{\alpha_3+1})
	\|\hat v_j-\hat v_i\|_{\alpha_1+1}^{\alpha_1+1}=0,
	\quad i,j\in[N].
	\]
	Because the trajectory is precompact, every sequence $t_n\to\infty$
	has a subsequence converging to the positive limit set. The continuity
	of the above quantities therefore gives
	\begin{align}
		\lim_{t\to\infty}\hat v_i(t)&=0,
		\qquad i\in\mathcal N_1,\label{C0-v}\\
		\lim_{t\to\infty}
		\kappa_f\phi_f(\|\hat x_i-\hat x_j\|_{\alpha_3+1})
		\|\hat v_j-\hat v_i\|_{\alpha_1+1}^{\alpha_1+1}&=0,
		\qquad i,j\in[N].\label{C0-diss}
	\end{align}
	We also record explicitly that the flocking-force terms themselves
	converge to zero. Set
	\[
	A_{ij}(t):=\kappa_f\phi_f(\|\hat x_i(t)-\hat x_j(t)\|_{\alpha_3+1}),
	\qquad
	r_{ij}(t):=\|\hat v_j(t)-\hat v_i(t)\|_{\alpha_1+1}.
	\]
	Both $A_{ij}$ and $r_{ij}$ are uniformly bounded by
	Proposition~\ref{P3.1} and the continuity of $\phi_f$. By norm
	equivalence,
	\[
	\|\overrightarrow{\mathrm{sig}}_{\alpha_1}
	(\hat v_j-\hat v_i)\|
	\leq C r_{ij}^{\alpha_1}.
	\]
	Moreover,
	\[
	A_{ij}r_{ij}^{\alpha_1}
	=
	\bigl(A_{ij}r_{ij}^{\alpha_1+1}\bigr)^{\frac{\alpha_1}{\alpha_1+1}}
	A_{ij}^{\frac1{\alpha_1+1}}.
	\]
	Thus~\eqref{C0-diss} and the boundedness of $A_{ij}$ imply
	\[
	\lim_{t\to\infty}
	\kappa_f\phi_f(\|\hat x_i-\hat x_j\|_{\alpha_3+1})
	\overrightarrow{\mathrm{sig}}_{\alpha_1}(\hat v_j-\hat v_i)
	=0,
	\qquad i,j\in[N].
	\]

	\noindent $\bullet$~(Proof of $\lim_{t \to \infty}\hat{x}_i(t)=0$):
	It follows from the above arguments that every element of the positive
	limit set {{$\omega(\hat{Z}(\cdot))$}} is of the following form:
	\begin{align*}
		\begin{aligned}
			&(\hat{x}_{i_1},\ldots,\hat{x}_{i_{N_1}},0,\ldots,0,\hat{x}_{j_1},\ldots,\hat{x}_{j_{N_2}},\hat{v}_{j_1},\ldots,\hat{v}_{j_{N_2}})\\
			&\text{satisfying}\quad
			\kappa_f\phi_f(\|\hat{x}_i-\hat{x}_j\|_{\alpha_3+1})
			\overrightarrow{\mathrm{sig}}_{\alpha_1}(\hat v_j-\hat v_i)=0,
			\quad i,j\in[N].
		\end{aligned}
	\end{align*}
	{{We now justify the corresponding acceleration limit.  For any
	$Y\in\omega(\hat Z(\cdot))$, semi-invariance
	(Definition~\ref{D2.3}) provides a solution through $Y$ that remains in
	$\omega(\hat Z(\cdot))$.  Along this solution, all successful-group
	velocities are identically zero, and hence their derivatives are zero as
	well.  Thus, the $\hat v_i$-component of the relative vector field
	vanishes at every $Y\in\omega(\hat Z(\cdot))$, $i\in\mathcal N_1$.
	Since the original relative trajectory is precompact and the vector
	field is continuous, any sequence $t_n\to\infty$ has a subsequence for
	which $\hat Z(t_n)$ converges to a point of the positive limit set.
	Consequently,}}
	\begin{align}\label{C1}
		{{\lim_{t \to \infty} \dot{\hat{v}}_i(t)  = 0, \quad \forall~ i\in \mathcal{N}_1.}}
	\end{align}
	Using the flocking-force limit established above, the convergence
	$\hat v_i(t)\to0$, and~\eqref{C1} in the second equation of
	\eqref{Fw-r}, we obtain, for every $i\in\mathcal N_1$,
	\[
	\lim_{t\to\infty}
	\psi_{r,i}(\|\hat x_i(t)\|_{\alpha_3+1}^{\alpha_3+1})
	\overrightarrow{\mathrm{sig}}_{\alpha_3}(\hat x_i(t))=0.
	\]
	Since $\hat x_i$ is uniformly bounded, norm equivalence gives a constant
	$C_i>0$ such that
	\[
	\|\hat x_i\|_{\alpha_3+1}^{\alpha_3+1}
	\leq C_i\|\overrightarrow{\mathrm{sig}}_{\alpha_3}(\hat x_i)\|
	\]
	along the trajectory. Consequently,
	\begin{align}\label{C2}
	\lim_{t\to\infty}
	\psi_{r,i}(\|\hat x_i(t)\|_{\alpha_3+1}^{\alpha_3+1})
	\|\hat x_i(t)\|_{\alpha_3+1}^{\alpha_3+1}=0.
	\end{align}
	We finally show that $\hat x_i(t)\to0$. Otherwise, for some
	$\varepsilon>0$ there would exist $t_n\to\infty$ such that
	$\|\hat x_i(t_n)\|_{\alpha_3+1}\geq\varepsilon$. By
	Proposition~\ref{P3.1}, the corresponding scalar arguments lie in the
	compact interval
	$[\varepsilon^{\alpha_3+1},U_i^{\alpha_3+1}]$. The continuity and strict
	positivity of $\psi_{r,i}$ on this interval imply
	\[
	\min_{s\in[\varepsilon^{\alpha_3+1},U_i^{\alpha_3+1}]}
	\psi_{r,i}(s)>0,
	\]
	which contradicts~\eqref{C2}. Hence
	\[
	\lim_{t\to\infty}\hat x_i(t)=0,
	\qquad i\in\mathcal N_1.
	\]
	{
		Since the choice of $Z(\cdot)\in \mathscr S(Z^0)$ is arbitrary, we conclude the desired theorem.
	}
}
\end{proof}

\subsection{Finite-time tracking success}\label{sec:3.3}
From Theorem \ref{T3.2}, we confirmed that under $(\mathcal{G}1)-(\mathcal{G}3)$, the pursuer configuration $\{(x_i, v_i)\}_{i \in \mathcal{N}_1}$ achieves asymptotic tracking success for the evader configuration $\{(x^*_i, v^*_i)\}_{i \in \mathcal{N}_1}$.
However, $(\mathcal{G}1)-(\mathcal{G}3)$ alone do not directly guarantee that $\{(x_i, v_i)\}_{i \in \mathcal{N}_1}$ can track the evader configuration $\{(x^*_i, v^*_i)\}_{i \in \mathcal{N}_1}$ within a finite time.
To address this, we consider an additional condition $(\mathcal{G}4)$ to $(\mathcal{G}1)$–$(\mathcal{G}3)$ and utilize the finite-time stability frameworks introduced in Section~\ref{sec:2.1} to obtain the desired result in this subsection.

Indeed, the finite-time tracking behavior observed in this model originates from the
non-Lipschitz structure of the coefficients in the governing equations.
This non-Lipschitz property is a key mechanism that enables finite-time
convergence.
In this context, the theorem below shows that, for any given initial data,
the pursuer configuration successfully tracks the evader configuration
for all $i \in \mathcal N_1$ within a fixed finite time that depends on
the initial data.
This property provides the direct motivation for the construction of
condition $(\mathcal{G}4)$.
Specifically, by prescribing the values of $\phi_f$ and $\psi_{r,i}$
over suitable intervals, condition $(\mathcal{G}4)$ is designed to ensure
that the subsystem of \eqref{Fw-r} associated with the pursuers that
successfully track the evader, together with the corresponding evaders,
exhibits finite-time stability.

{
	Before proving finite-time tracking success, we first identify the
	closed subsystem associated with the agents in $\mathcal{N}_1$ and
	derive the Lyapunov stability of its origin.
	By Proposition~\ref{P3.1}, Theorem~\ref{T3.1}, the strict
	separation and support conditions in $(\mathcal{G}4)$, and the fifth assertion of Remark \ref{R3.1}, one has
	\[
	\phi_f(\|\hat{x}_i-\hat{x}_j\|_{\alpha_3+1})=0,
	\quad
	i\in\mathcal{N}_1,
	\quad
	j\in\mathcal{N}_2.
	\]
	Therefore, the dynamics of
	$\{(\hat{x}_i,\hat{v}_i)\}_{i\in\mathcal{N}_1}$ form the following
	closed subsystem:
	\begin{equation}\label{Fw-r-N1}
		\begin{cases}
			\displaystyle
			\dot{\hat{x}}_i=\hat{v}_i,
			\quad t>0,
			\quad i\in\mathcal{N}_1,
			\\[0.2cm]
			\displaystyle
			\dot{\hat{v}}_i
			=
			\frac{\kappa_f}{N}
			\sum_{j\in\mathcal{N}_1}
			\phi_f(\|\hat{x}_i-\hat{x}_j\|_{\alpha_3+1})
			\overrightarrow{\textup{sig}}_{\alpha_1}
			(\hat{v}_j-\hat{v}_i)
			\\[0.2cm]
			\displaystyle\hspace{1cm}
			-(\gamma_i-\gamma_i^*)
			\overrightarrow{\textup{sig}}_{\alpha_2}(\hat{v}_i)
			-
			\psi_{r,i}
			(\|\hat{x}_i\|_{\alpha_3+1}^{\alpha_3+1})
			\overrightarrow{\textup{sig}}_{\alpha_3}(\hat{x}_i),
			\\[0.2cm]
			\displaystyle
			(\hat{x}_i,\hat{v}_i)\big|_{t=0+}
			=
			(\hat{x}_i(0),\hat{v}_i(0))
			\in\mathbb{R}^{2d}.
		\end{cases}
	\end{equation}
	
	We now study the Lyapunov stability of the origin for the closed
	subsystem~\eqref{Fw-r-N1}.
	
{
	\begin{lemma}\label{L3.2}
		\emph{(Lyapunov stability of the successful-tracking subsystem)}
		Suppose that the framework $(\mathcal{G}1)-(\mathcal{G}4)$ holds.
		Then, the origin of the closed subsystem~\eqref{Fw-r-N1} is
		Lyapunov stable. More precisely, for every $\varepsilon>0$, there
		exists $\delta>0$ such that every {{maximal}} classical solution
		\[
		\widehat{Z}_1(\cdot)
		=
		\{(\hat{x}_i(t),\hat{v}_i(t))\}_{i\in\mathcal{N}_1}
		\]
		to the closed subsystem~\eqref{Fw-r-N1}, {{initially defined on
		$[0,T_{\max})$,}} and satisfying
		\[
		\left(
		\sum_{i\in\mathcal{N}_1}
		\left(
		\|\hat{x}_i(0)\|_{\alpha_3+1}^{2}
		+
		\|\hat{v}_i(0)\|^{2}
		\right)
		\right)^{\frac12}
		<\delta
		\]
		{{is global}} and satisfies
		\[
		\left(
		\sum_{i\in\mathcal{N}_1}
		\left(
		\|\hat{x}_i(t)\|_{\alpha_3+1}^{2}
		+
		\|\hat{v}_i(t)\|^{2}
		\right)
		\right)^{\frac12}
		<\varepsilon,
		\quad t\geq0.
		\]
	\end{lemma}
	
	\begin{proof}
		Let
		\[
		\widehat{Z}_1(\cdot)
		=
		\{(\hat{x}_i(t),\hat{v}_i(t))\}_{i\in\mathcal{N}_1}
		\]
		be any {{maximal}} classical solution to the closed
		subsystem~\eqref{Fw-r-N1} {{on $[0,T_{\max})$}}. We define the associated
		energy functional by
		\[
		\mathcal{E}_1(t)
		:=
		\frac12\sum_{i\in\mathcal{N}_1}\|\hat{v}_i(t)\|^2
		+
		\frac{1}{\alpha_3+1}
		\sum_{i\in\mathcal{N}_1}
		\int_0^{\|\hat{x}_i(t)\|_{\alpha_3+1}^{\alpha_3+1}}
		\psi_{r,i}(u)\,du.
		\]
		{{By the same calculation as in Lemma~\ref{L3.1}, throughout the
		maximal interval one has}}
		\[
		\frac{d\mathcal{E}_1}{dt}
		=
		-\frac{\kappa_f}{2N}
		\sum_{i,j\in\mathcal{N}_1}
		\phi_f(\|\hat{x}_i-\hat{x}_j\|_{\alpha_3+1})
		\|\hat{v}_j-\hat{v}_i\|_{\alpha_1+1}^{\alpha_1+1}
		-
		\sum_{i\in\mathcal{N}_1}
		(\gamma_i-\gamma_i^*)
		\|\hat{v}_i\|_{\alpha_2+1}^{\alpha_2+1}
		\leq0.
		\]
		Moreover, since
		\[
		\psi_{r,i}(0)\geq0,
		\quad
		\psi_{r,i}(u)>0
		\quad\text{for all }u>0,
		\]
		the functional $\mathcal{E}_1$ is continuous and positive definite
		with respect to the origin. Next, let
		\[
		\widehat{Z}_1
		:=
		\left(
		\{\hat{x}_i\}_{i\in\mathcal{N}_1},
		\{\hat{v}_i\}_{i\in\mathcal{N}_1}
		\right),
		\]
		and set
		\[
		\|\widehat{Z}_1\|_*
		:=
		\left(
		\sum_{i\in\mathcal{N}_1}
		\left(
		\|\hat{x}_i\|_{\alpha_3+1}^{2}
		+
		\|\hat{v}_i\|^{2}
		\right)
		\right)^{\frac12}.
		\]
		For a given $\varepsilon>0$, {{define}}
		\[
		m_\varepsilon
		:=
		\min_{\|\widehat{Z}_1\|_*=\varepsilon}
		\mathcal{E}_1(\widehat{Z}_1).
		\]
		The continuity and positive definiteness of $\mathcal{E}_1$ imply
		that $m_\varepsilon>0$. {{Choose $\delta\in(0,\varepsilon)$
		sufficiently small that}}
		\[
		\|\widehat{Z}_1(0)\|_*<\delta
		\quad\Longrightarrow\quad
		\mathcal{E}_1(0)<m_\varepsilon.
		\]
		Since
		\[
		\mathcal{E}_1(t)\leq\mathcal{E}_1(0),
		{{\quad 0\leq t<T_{\max},}}
		\]
		the solution cannot reach the {{sphere}}
		$\{\|\widehat{Z}_1\|_*=\varepsilon\}$. Consequently,
		\[
		\|\widehat{Z}_1(t)\|_*<\varepsilon,
		\quad 0\leq t<T_{\max}.
		\]
		{{Hence the solution remains in a compact subset of the state space on
		its maximal interval. Since the vector field of~\eqref{Fw-r-N1} is
		continuous, the continuation theorem gives $T_{\max}=\infty$.
		Since the maximal classical solution was arbitrary, the desired lemma
		follows.}}
	\end{proof}
}
}

{
\begin{lemma}\label{L3.3}
\emph{(Uniform finite-time stability on the initial-energy sublevel set)}
Fix an admissible initial datum $Z^0$ and set
\[
c_*:=\mathcal E(0).
\]
If $c_*>0$, define
\[
K_*
:=
\left\{
\widehat Z_1:\mathcal E_1(\widehat Z_1)\leq c_*
\right\},
\]
where $\mathcal E_1$ is the successful-subsystem energy introduced in
Lemma~\ref{L3.2}. Under $(\mathcal G1)-(\mathcal G4)$, there exists a
finite constant $T_*>0$, depending only on $Z^0$ and the parameters, such
that every maximal classical solution of~\eqref{Fw-r-N1} with initial
value in $K_*$ is global and satisfies
\[
\widehat Z_1(t)=0,
\qquad t\geq T_*.
\]
The same $T_*$ applies to all possibly non-unique solutions.
\end{lemma}

\begin{proof}
We verify the hypotheses of Lemma~\ref{L2.1}. First, $K_*$ is closed and
contains the origin in its interior. Since the terms in $\mathcal E_1$ are
nonnegative and $c_*=\mathcal E(0)$, the definition of $U_i$ in
$(\mathcal G2)$ gives
\[
\|\hat x_i\|_{\alpha_3+1}\leq U_i,
\qquad i\in\mathcal N_1,
\qquad \widehat Z_1\in K_*.
\]
The kinetic part bounds all $\hat v_i$, so $K_*$ is compact. The energy
identity in Lemma~\ref{L3.2} is valid on every maximal interval of
existence and gives
\[
\frac{d\mathcal E_1}{dt}\leq0.
\]
Consequently, every maximal solution starting in $K_*$ remains in $K_*$
throughout its maximal interval. Since the vector field is continuous and
bounded on the compact set $K_*$, the continuation theorem makes every
such maximal solution global. Thus $K_*$ is forward positively invariant
in the sense of Definition~\ref{D2.1}.

We next prove asymptotic stability relative to $K_*$. Lyapunov stability
follows from Lemma~\ref{L3.2}. Let $\widehat Z_1(\cdot)$ be any solution
starting in $K_*$. Its trajectory is precompact, and $\mathcal E_1$ is
locally Lipschitz. The nonsmooth LaSalle principle therefore yields
\[
\omega(\widehat Z_1(\cdot))
\subset
\left\{
\widehat Z_1:\hat v_i=0\text{ for all }i\in\mathcal N_1
\right\},
\]
because $\gamma_i-\gamma_i^*>0$ on $\mathcal N_1$. If
$Y\in\omega(\widehat Z_1(\cdot))$, semi-invariance provides a solution
through $Y$ that remains in this positive limit set. Along that solution,
all $\hat v_i$ vanish identically, and hence so do their derivatives. The
second equation of~\eqref{Fw-r-N1} then gives
\[
\psi_{r,i}(\|\hat x_i\|_{\alpha_3+1}^{\alpha_3+1})
\overrightarrow{\mathrm{sig}}_{\alpha_3}(\hat x_i)=0,
\qquad i\in\mathcal N_1.
\]
By $(\mathcal G2)$, $\psi_{r,i}(u)>0$ for $u>0$, so $\hat x_i=0$ for every
$i\in\mathcal N_1$. Thus the positive limit set is $\{0\}$, and
precompactness implies $\widehat Z_1(t)\to0$. The origin is therefore
asymptotically stable relative to $K_*$. 

It remains to verify local negative homogeneity. For the alternative of
$(\mathcal G4)$ under consideration, define
\[
q_f:=
\begin{cases}
0,&\text{in }(\mathcal G4)(1)\text{ or }(\mathcal G4)(2),\\
\sigma,&\text{in }(\mathcal G4)(3)\text{ or }(\mathcal G4)(4),
\end{cases}
\qquad
q_r:=
\begin{cases}
0,&\text{in }(\mathcal G4)(1)\text{ or }(\mathcal G4)(3),\\
\beta,&\text{in }(\mathcal G4)(2)\text{ or }(\mathcal G4)(4).
\end{cases}
\]
Then the exponent relations in $(\mathcal G4)$ are
\[
b=a+\lambda,
\qquad
b\alpha_1+aq_f=b\alpha_2=b+\lambda
=aq_r(\alpha_3+1)+a\alpha_3,
\qquad \lambda<0.
\]
Consider
\[
\delta_\varepsilon(\widehat Z_1)
:=
\left(
\{\varepsilon^a\hat x_i\}_{i\in\mathcal N_1},
\{\varepsilon^b\hat v_i\}_{i\in\mathcal N_1}
\right),
\qquad \varepsilon\in(0,1].
\]
For every $\widehat Z_1\in K_*$,
\[
\|\hat x_i-\hat x_j\|_{\alpha_3+1}
\leq2\max_{\ell\in\mathcal N_1}U_\ell,
\qquad
\|\hat x_i\|_{\alpha_3+1}^{\alpha_3+1}
\leq U_i^{\alpha_3+1}.
\]
Hence the local scaling identities in $(\mathcal G4)$ are valid on $K_*$
and under every contractive dilation. They give
\[
\phi_f(\varepsilon^a u)=\varepsilon^{aq_f}\phi_f(u),
\qquad
\psi_{r,i}(\varepsilon^{a(\alpha_3+1)}s)
=\varepsilon^{a(\alpha_3+1)q_r}\psi_{r,i}(s).
\]
Together with
$\overrightarrow{\mathrm{sig}}_\alpha(\varepsilon z)
=\varepsilon^\alpha\overrightarrow{\mathrm{sig}}_\alpha(z)$, these
relations show that the vector field of~\eqref{Fw-r-N1} is homogeneous of
degree $\lambda<0$ on $K_*$. Therefore Lemma~\ref{L2.1}, applied with
$K=K_*$, yields a neighborhood $V\subset K_*$ of the origin and a finite
local settling-time bound $T_{\mathrm{loc}}>0$ that is valid for every
solution starting in $V$.

It remains only to make the entrance into $V$ uniform over all solutions
starting in $K_*$. We claim that there exists $T_{\mathrm{ent}}>0$ such
that every solution starting in $K_*$ enters $V$ no later than
$T_{\mathrm{ent}}$. If not, for each $n\in\mathbb N$ there would be a
solution $\widehat Z_1^n(\cdot)$ starting in $K_*$ and satisfying
\[
\widehat Z_1^n(t)\notin V,
\qquad 0\leq t\leq n.
\]
All these trajectories remain in the compact invariant set $K_*$, and
the vector field is bounded there. Arzel\`a--Ascoli and a diagonal
argument therefore provide a subsequence converging uniformly on compact
time intervals to a global solution $\widehat Z_1^\infty(\cdot)$ starting
in $K_*$. Since the origin is attractive relative to $K_*$, this limit
solution belongs to $V$ at some finite time. The openness of $V$ and
uniform convergence then imply that $\widehat Z_1^n$ also enters $V$
before time $n$ for all sufficiently large $n$, a contradiction.

Once a given solution enters $V$, its time shift is a global solution
starting in $V$, so Lemma~\ref{L2.1} makes it identically zero after at
most an additional time $T_{\mathrm{loc}}$. Thus
\[
T_*:=T_{\mathrm{ent}}+T_{\mathrm{loc}}
\]
is a common settling-time bound for all initial values in $K_*$ and all
possibly non-unique solutions. This proves the lemma.
\end{proof}
}

{
We are now ready to prove finite-time tracking success with a settling
time that is common to all solutions issued from the same initial datum.

\begin{theorem}\label{T3.3}~\emph{(Finite-time tracking success)}
For a given admissible initial datum $Z^0$, suppose that
$(\mathcal G1)-(\mathcal G4)$ hold.  Then there exists a finite time
$T=T(Z^0)\in[0,\infty)$ such that every global classical solution
\[
Z(\cdot)
=
\{(x_i(t),v_i(t),x_i^*(t),v_i^*(t))\}_{i=1}^N
\in\mathscr S(Z^0)
\]
satisfies
\[
x_i(t)-x_i^*(t)=0,
\qquad
v_i(t)-v_i^*(t)=0,
\qquad
i\in\mathcal N_1,\quad t\geq T(Z^0).
\]
\end{theorem}

\begin{proof}
The strict separation and support assumptions in $(\mathcal G4)$,
together with Proposition~\ref{P3.1} and Theorem~\ref{T3.1}, imply that
for every global solution issued from $Z^0$,
\[
\phi_f(\|\hat x_i-\hat x_j\|_{\alpha_3+1})=0,
\qquad i\in\mathcal N_1,\quad j\in\mathcal N_2.
\]
Thus the successful-group variables solve the same closed
subsystem~\eqref{Fw-r-N1}, with the common initial state
$\widehat Z_1^0$ determined by $Z^0$.

Set $c_*:=\mathcal E(0)$. If $c_*=0$, then the nonnegativity of all
terms in the total energy gives $\mathcal E_1(\widehat Z_1^0)=0$.
Positive definiteness and monotonicity of $\mathcal E_1$ imply
$\widehat Z_1(t)=0$ for every solution and all $t\geq0$, so we may take
$T(Z^0)=0$. If $c_*>0$, then
\[
\mathcal E_1(\widehat Z_1^0)\leq\mathcal E(0)=c_*,
\]
so $\widehat Z_1^0\in K_*$. Lemma~\ref{L3.3} provides a finite constant
$T_*$ that is valid for every solution of the closed subsystem starting
in $K_*$. In particular, the projection of every full solution in
$\mathscr S(Z^0)$ satisfies
\[
\hat x_i(t)=\hat v_i(t)=0,
\qquad i\in\mathcal N_1,\quad t\geq T_*.
\]
Taking $T(Z^0):=T_*$ proves the theorem. The set $K_*$ and the resulting
bound depend only on the fixed initial datum $Z^0$ and the parameters, not
on the particular solution selected from $\mathscr S(Z^0)$.
\end{proof}
}

\section{Conclusion} \label{sec:4}
\setcounter{equation}{0}
In this study, we have proposed a sufficient framework to determine the success or failure of finite-time tracking in pursuer--evader dynamics.
To this end, we first designed a nonsmooth coupled multi-agent system in which the pursuers try to track the evaders within a finite time, the evaders try to avoid collisions with the pursuers, and both groups move in an organized manner under suitable control inputs.
The proposed system was then analyzed by using basic results on finite-time stability for nonsmooth systems, several energy estimates, and the nonsmooth version of LaSalle’s invariance principle.
As a result, we established sufficient conditions, in terms of initial data, system parameters, and kernel functions, that guarantee either finite-time tracking success or failure between the pursuers and evaders.
This work provides a mathematical framework for systematically analyzing the success or failure of finite-time tracking in nonsmooth multi-agent systems with complex interactions.
However, the present results focus on theoretical analysis for a finite number of agents in Euclidean space.
An important direction for future research is to extend this framework to the continuum or mean-field setting, or to models defined on Riemannian manifolds. Accordingly, we plan to extend the nonsmooth pursuer--evader model proposed in this paper either to a kinetic model or to a generalized formulation on Riemannian manifolds, in order to explore richer geometric and analytical properties.

\section*{Declarations}

\noindent\textbf{Conflict of Interest}

\medskip

The authors declare that they have no conflict of interest.

\vspace{0.5cm}

\noindent\textbf{Data Availability}

\medskip

No data were generated or analyzed in this study.

\end{document}